\documentclass[sn-mathphys-num]{sn-jnl}
\usepackage[utf8]{inputenc}
\usepackage{a4wide}
\usepackage{amssymb}
\usepackage{amsmath}
\usepackage{amsthm}
\usepackage{hyperref}
\usepackage{breqn}
\usepackage{xcolor}
\usepackage{mathtools}
\usepackage[numbers]{natbib}
\usepackage{bm}
\usepackage{enumitem}
\usepackage{mathrsfs}
\usepackage{graphicx}
\usepackage{caption}
\usepackage{url}
\usepackage{placeins}
\usepackage[english]{babel}
\usepackage{array}
\usepackage{booktabs}
\usepackage{subcaption}
\usepackage{multirow}
\numberwithin{equation}{section}
\newcommand{\dm}[1]{{\displaystyle{#1}}}
\newtheorem{corollary}{Corollary}[section]
\newtheorem{lemma}{Lemma}[section]
\theoremstyle{remark}
\newtheorem{exam}{\bf Example}[section]
\def \vec{\mathrm v\mathrm e \mathrm c}
\def \R{{\mathbb R}}

\def \x{{\bm x}}
\def \b{{\bf b}}

\def \T{\mathsf T}
\def\bmatrix#1{\left[\begin{matrix}
		#1
	\end{matrix}\right]}

\def \diag{\mathrm{diag}}

\def \D{{\Delta}}

\def \be{\bm{\chi}}

\def  \0{\bf 0}

\def \rank{\mathrm{rank}}

\def \x{\bm x}
\def \l{\bm l}
\def \y{\bm y}
\def \z{\bm z}
\def \T{{\bm\top}}
\def \vp{\bm \varphi}

\def\bmatrix#1{\left[ \begin{matrix} #1 \end{matrix} \right]}

\def \R{{\mathbb R}}

\theoremstyle{thmstyleone}%
\newtheorem{theorem}{Theorem}
\theoremstyle{thmstyletwo}%
\newtheorem{remark}{Remark}%

\theoremstyle{thmstylethree}%
\newtheorem{definition}{Definition}%

\begin{document}

\title[Partial Condition Numbers for Double Saddle Point Problems]{Partial Condition Numbers for Double Saddle Point Problems}


\author*[1]{\fnm{Sk. Safique} \sur{Ahmad}}\email{safique@iiti.ac.in}

\author[1]{\fnm{Pinki} \sur{Khatun}}\email{pinki996.pk@gmail.com}


\affil*[1]{\orgdiv{Department of Mathematics}, \orgname{Indian Institute of Technology Indore}, \orgaddress{\street{Simrol}, \city{Indore}, \postcode{453552}, \state{Madhya Pradesh}, \country{India}}}

\abstract{This paper presents a unified framework for investigating the partial condition number (CN) of the solution of double saddle point problems (DSPPs) and provides closed-form expressions for it. This unified framework encompasses the well-known partial normwise CN (NCN), partial mixed CN (MCN) and partial componentwise CN (CCN) as special cases. Furthermore, we derive sharp upper bounds for the partial NCN, MCN and CCN, which are computationally efficient and free of expensive Kronecker products. By applying perturbations that preserve the structure of the block matrices of the DSPPs, we analyze the structured partial NCN, MCN and CCN  when the block matrices exhibit linear structures. By leveraging the relationship between DSPP and equality constrained indefinite least squares (EILS) problems, we recover the partial CNs for the EILS problem. Numerical results confirm the sharpness of the derived upper bounds and demonstrate their effectiveness in estimating the partial CNs. }





\maketitle

\section{Introduction}
We consider the following linear system  with the double saddle point structure:
\begin{align}\label{eq1:SPP}
   \mathfrak{B}\bm{w}=\bf{b},~ 
\end{align}
 where  
 \begin{eqnarray}
     \mathfrak{B}=\bmatrix{A & B^{\bm\top} & \bf 0\\ B &-D & C^{\bm\top} \\ \bf 0 & C& E}; ~ \bm{w}=\bmatrix{\bm{x}\\ \bm{y} \\ \bm{z}};~ \bf{b}=\bmatrix{\mathbf{b}_1\\\mathbf{b}_2 \\ \mathbf{b}_3};
 \end{eqnarray}
 $A\in \R^{n\times n},$ $B\in \R^{m\times n},$ $C\in \R^{p\times m},$ $D\in \R^{m\times m},$ $E\in \R^{p\times p},$  {$\mathbf{b}\in \R^{\l}$} and $\l=n+m+p.$ Then the linear system of the form \eqref{eq1:SPP} is known as the double saddle point problem (DSPP). We refer  $\mathfrak{B}$ as the double saddle point matrix. Conditions on the invetability for the matrix $\mathfrak{B}$ have been studied in \cite{beik2024, Greif2023}. 
 To ensure a unique solution to \eqref{eq1:SPP}, throughout the paper, we assume that $\mathfrak{B}$ is nonsingular.

Linear systems of the form \eqref{eq1:SPP} emerge in numerous computational science and technological applications, such as in finite element discretizations of coupled Stokes-Darcy flow equations \cite{stokesdarcy2019}, the equality constrained indefinite least square problem (EILS)  \cite{ILSE2003}, computational fluid dynamics \cite{CFD2005},  liquid crystal director modeling \cite{LCDM2014}, PDE-constraint optimization problems \cite{PDE-constrained2010}, and so on. Due to the versatile applicability of the DSPP \eqref{eq1:SPP}, their numerical solutions have garnered significant interest and have been extensively investigated in recent years. Several iterative methods and preconditioning techniques have been developed in recent times for solving DSPPs; see \cite{ Greif2023, HuangNA, Multiparameter2023, Pinki_PESS} and references therein.

 The coefficient matrix of a DSPP is typically obtained from measurements and is therefore subject to some measurement errors, as well as subject to the rounding and truncation errors introduced due to the finite precision of the computers. As a result of these inevitable perturbations, we work with a matrix that is slightly different from the original one. Since such errors are often unavoidable, it is important to investigate how small perturbations in the input matrices influence the resulting solution.
 Perturbation analysis is one of the most useful tool in numerical analysis for assessing the sensitivity and stability of the solution of a problem obtained through numerical methods \cite{Higham2002}. The condition number (CN) quantifies the worst-case sensitivity of a solution to a small perturbation in the input data matrices. Additionally, CNs are vital in providing first-order error estimates for the relative error of a numerical solution.  The general framework of CNs was first introduced by \citet{Rice1966}. It primarily focuses on the normwise CN (NCN), where norms are used to quantify both the input perturbations and the resulting errors in the output data. However, a significant limitation of NCN is its inability to take into account the inherent structure of poorly scaled or sparse input data. To overcome this limitation, the mixed CN (MCN) and componentwise CN (CCN) have gained increasing attention in the literature \cite{Gohberg1993, rohn1989new, skeel1979scaling}. The former quantifies input perturbations componentwise and output errors using norms, while the later measures both input perturbations and output errors componentwise.

 {Perturbation analysis, CNs, and structured backward error analysis for standard saddle point problems (SPPs) of the following form:
\begin{align}\label{GSPP}
\bmatrix{
F & G^{\T}\\ 
H & K}
\bmatrix{
u\\ 
p}
=
\bmatrix{
f\\ 
g
}   
\end{align}
have been extensively studied in the literature.
 For example, the NCN 
for the solution of the SPP \eqref{GSPP}, under the assumptions \( F = F^{\T} \), \( G = H \), and \( K = \0 \), was initially investigated in~\cite{conditionwang}. In~\cite{Xiangcnd2007}, explicit expressions for the NCN, MCN, and CCN were derived for the solution vector $[u^{\T},\,p^{\T}]^{\T}$, as well as for its components \( u \) and \( p \), under the same structural assumptions. Building on this,~\citet{weiwei2009} further examined the NCN and established perturbation bounds for the combined solution vector \( [u^\T, \, p^\T]^\T \). Subsequently, \citet{meng2019condition} extended this line of research by analyzing both the MCN and CCN for the full solution and its individual components \( u \) and \( p \). More recently, structured CNs and backward error analysis for the SPP~\eqref{GSPP} have been explored in~\cite{PinkiGSPP, Pinki_BE}, particularly when the coefficient matrices exhibit linear structures such as Toeplitz or symmetry. However, these studies fail to leverage the special three-by-three block structure of the coefficient matrix \( \mathfrak{B} \), and the CNs for the DSPP \eqref{eq1:SPP} remain unexplored and still require further detailed investigation.
 }

In recent years, the structured CNs of various problems have been studied, emphasizing the preservation of the linear structure of the original matrices in the perturbation matrices; see, for example \cite{Rump2003, RumpII, cucker2007mixed2, PinkiGSPP}. 
 The block matrices \(A\), \(D\) and \(E\) often exhibit particular linear structures in various applications; see \cite{PDE-constrained2010, LCDM2014}. This makes it compelling to explore the structured CNs for the DSPP \eqref{eq1:SPP} by preserving the linear structures of the diagonal block matrices to their corresponding perturbation matrices.


In many applications, the components \(\x\), \(\y\), and \(\z\) of the solution vector \(\bm{w}\) correspond to distinct physical quantities. For example, in PDE-constrained optimization problems \cite{PDE-constrained2010}, \(\x\) typically represents the desired state, \(\z\) the control variable, and \(\y\) the Lagrange multipliers. Since each component may respond differently to perturbations, it is essential to analyze the sensitivity of individual components rather than the full solution vector \cite{Cao2003}. Traditional CNs, which capture the overall sensitivity, lack the ability to reveal the conditioning of a specific part of the solution. To overcome this limitation, the concept of a partial or projected CN was introduced, which quantifies the sensitivity of a linear function of the solution using a selection matrix \(\mathrm{L}\). This matrix allows one to isolate or select and study the sensitivity of a specific solution component. First proposed in \cite{Cao2003} for linear systems, this framework has since been widely extended to various problems, including linear least squares \cite{LLS2007}, weighted least squares \cite{WLS2017}, indefinite least squares \cite{ILS2018, ILSE2022}, total least squares \cite{TLS2011}, and generalized SPPs \cite{PinkiGSPP}.

In this paper, we consider the CN of the linear function $\mathrm{L}[\x^{\bm\top},\y^{\bm\top},\z^{\bm\top}]^{\bm\top}$ of the solution $\bm{w}=[\x^{\bm\top},\y^{\bm\top},\z^{\bm\top}]^{\bm\top},$ where $\mathrm{L}\in \R^{k\times {\bm l}}$ $(k\leq {\bm l})$ or the partial CN of the solution $\bm{w}=[\x^{\bm\top},\y^{\bm\top},\z^{\bm\top}]^{\bm\top}$ of the DSPP \eqref{eq1:SPP}. Different choices of \( \mathrm{L} \) provide the flexibility to determine the CNs of various solution components of \( \bm{w} \). For example, by selecting \(\mathrm{L} = I_{\bm{l}}\), \([I_n~~ \mathbf{0}_{n \times (\bm{l}-n)}]\), or $[{\bf 0}~~I_m~~ \mathbf{0}_{m \times p}]$, we can determine the CN of \(\bm{w}\), \(\x\), or  \(\bm{y}\), respectively. Furthermore, our investigation presents a general framework that encompasses well-known CNs, such as NCN, MCN and CCN, as special cases.

The key contributions of the paper are highlighted as follows:
\begin{itemize}
    \item In this work, we explore a general form of partial CN, referred to as unified partial CN, which has a versatile nature and provides a comprehensive framework encompassing the partial NCN, MCN, and CCN of the solution of the DSPP.
    \item By considering structure-preserving perturbations on $A,$ $D$, and $E,$ when they retain some linear structures, we derive structured partial CNs for the DSPP.
    \item By exploring the connection between DSPPs and EILS problems, we demonstrate that our derived CN formula can be used to recover the partial CNs for the EILS problem.
    \item Numerical experiments demonstrate that the derived upper bounds provide sharp estimates of the partial CNs. Furthermore, the partial CNs offer precise estimates of the relative forward error in the solution.
\end{itemize}

The structure of the rest of the paper is as follows. Section \ref{sec2}  introduces a few notations, basic definitions, and preliminaries. Section \ref{sec:PCN} presents a unified framework partial CN of the solution of the DSPP \eqref{eq1:SPP}. Section \ref{sec:SPCN} focuses on the investigation of the structured partial CNs for the DSPP. In Section \ref{sec:EILS}, we discuss the partial CNs for EILS problems. Section \ref{sec:Numerical}
 consists of some numerical examples. Section \ref{sec:Conclusion} includes the concluding statements.
 
 \section{Notation and preliminaries}\label{sec2}
\subsection{Notation}
Throughout the paper, $\R^{m\times n}$ represents the set of all $m\times n$ real matrices. The symbols $\|\cdot \|_F$, $\|\cdot \|_2,$ and $\|\cdot \|_{\infty}$ stand for the  Frobenius norm, spectral norm, and infinity norm of the argument, respectively. The notation ${\bf 1}_n\in \R^n$ represents a vector with all entries equal to $1$. For any vector $z=[z_1,z_2,\ldots, z_n]^{\bm\top}\in \R^n,$ $\Theta_{z}\in \R^{n\times n}$ is the diagonal matrix with $i$-th diagonal entry $z_i.$  Following \cite{polarD}, for any vector $z\in \R^n,$ we define 
\begin{align}
   z^{\ddagger}=[z^{\ddagger}_1,z^{\ddagger}_2,\ldots, z^{\ddagger}_n]^{\bm\top},
\end{align}
 where \begin{align}
     z_i^{\ddagger}=\left\{\begin{array}{cc}
   \frac{1}{z_i},      &  z_i\neq 0, \\
      1,    &   z_i=0.
     \end{array}\right.
 \end{align}
 Moreover, the entrywise division of two vector $z,w\in \R^n$ is defined as follows:
 \begin{equation*}
     \frac{z}{w}=\Theta_{w^{\ddagger}}z.
 \end{equation*}
 For the matrices $X=[x_{ij}]\in \R^{m\times n}$ and $Y=[y_{ij}]\in \R^{m\times n},$ we define $X\odot Y:=[x_{ij}y_{ij}]\in \R^{m\times n}.$ Note that, for $z,w\in \R^{n},$ $\dm{\frac{z}{w}}=w^{\ddagger}\odot z.$   For  $X \in \R^{m\times n},$ we set $|X|:=[|x_{ij}|].$   For $X,\,Y\in\R^{m\times n},$  the notation $|X|~
\leq |Y|$ represents $|x_{ij}|~\leq |y_{ij}|$ for all $1\leq i\leq m$ and $1\leq j\leq n$. For  $X=[\x_1, \x_2,\ldots, \x_n]\in \R^{m\times n},$ where $\x_i\in \R^{m},$ $i=1,2,\ldots, n,$  we define $\vec(X):=[\x_1^{T}, \x_2^{T},\ldots, \x_n^{T}]^{T}\in\R^{mn}$. For given matrices $A_1,A_2,\ldots, A_n$, we use $$\vec(\mathbf{X}):=[\vec(A_1)^{\bm\top},~ \vec(A_2)^{\bm\top}, ~\ldots,~ \vec(A_n)^{\bm\top}]^{\bm\top},$$ where $\mathbf{X}=(A_1,A_2,\ldots, A_n).$ The Kronecker product \cite{kronecker1981} of  $X\in \R^{m\times n}$ and $Y\in \R^{p\times q}$ is defined as $X\otimes Y:= [x_{ij}Y]\in \R^{mp\times nq}.$ For matrices $A, B, X, Y$ and $Z$ with compatible dimensions, we have the following properties \cite{kronecker1981, kronecker2004}:
 \begin{align}\label{kron}
\left\{ \begin{array}{c}
 \vec(AZB) =  (B^{\bm\top}\otimes A)\vec(Z),   \\
   (X\otimes Y)^{\bm\top}=X^{\bm\top}\otimes Y^{\bm\top},   \\
   |X\otimes Y|=|X|\otimes |Y|,\\
   (A\otimes B)(X\otimes Y)=AX\otimes BY.\\
 \end{array}\right.
  \end{align}

 \subsection{Preliminaries}
At the beginning of this subsection, we introduce the concept of the general CN, referred to as the unified CN. 
\begin{definition}\cite{UCN2018}\label{def:ucn}
    Let $\Upsilon : \R^p \mapsto \R^q$ be a continuous mapping defined on an open set $\Omega_{\Upsilon}\subseteq \R^p.$ Then, the unified CN of $\Upsilon$ at $\bm v\in \Omega_{\Upsilon}$ is defined by 
    \begin{equation}
        \mathfrak{K}_{\Upsilon}(\bm v) = \lim_{\bm{\epsilon} \rightarrow 0}\sup_{0< \|\chi^{\ddagger} \odot \D \bm v\|_{\tau}\leq \bm{\epsilon}} \frac{\left\|\xi^{\ddagger}\odot \left(\Upsilon(\bm v+\D \bm v)-\Upsilon(\bm v)\right)\right\|_{\gamma}}{\|\chi^{\ddagger} \odot \D \bm v\|_{\tau}},
    \end{equation}
where $\xi \in \R^{q},$ ${\chi}\in \R^p$ are the parameters such that if some entry of $\chi$ is zero, then the corresponding entry of $\D \bm v$ must be zero, and $\|\,\cdot\,\|_{\tau}$ and $\|\,\cdot\,\|_{\gamma}$ are two vector norms defined on $\R^p$ and $\R^q,$ respectively.
\end{definition}

The parameters $\chi$ and $\xi$ in Definition \ref{def:ucn} can also be positive real numbers, rather than vectors. Note that, Definition \ref{def:ucn} leads to  the following bound:
\begin{equation}\label{eq:forward}
\left\|\xi^{\ddagger}\odot (\Upsilon(\bm v+\D \bm v)-\Upsilon(\bm v))\right\|_{\gamma}\leq \mathfrak{K}_{\Upsilon}(\bm v) \|\chi^{\ddagger} \odot \D \bm v\|_{\tau}\,+~\mathcal{O}(\|\chi^{\ddagger} \odot \D \bm v\|^2_{\tau}).
\end{equation}
Therefore, the forward error in the solution can be estimated using CNs.
\begin{remark}\label{re21}
   The unified CN described in Definition \ref{def:ucn} represents a broad generalization of various well-known CNs that have been explored in the literature. For example:
   \begin{itemize}
       \item   \textbf{NCN:} Consider $\tau=\gamma=2,$
       $\chi =\|\bm v\|_2{\bf 1}_p$ $\in \R^p$ with $\bm v\neq \0$ and $\xi=\|\Upsilon(\bm v)\|_2{\bf 1}_q\in \R^q$ with $\Upsilon(\bm v)\neq \0,$ then  we obtain the NCN, denoted by $ \mathfrak{K}^{(2)}_{\Upsilon}(\bm v).$

 \item    \textbf{MCN:} Consider $\tau=\gamma=\infty,$ $\chi=\bm v\neq \0,$  $\xi=\|\Upsilon(\bm v)\|_{\infty}{\bf 1}_q\in \R^q$ with $\Upsilon(\bm v)\neq \0,$ then  we obtain the MCN, denoted by $ \mathfrak{K}^{\infty}_{mix,\Upsilon}(\bm v).$

\item    \textbf{CCN:} Consider $\tau=\gamma=\infty,$ $\chi =\bm v\neq \0,$  and $\xi=\Upsilon(\bm v)\in \R^q$ with  $\Upsilon(\bm v)\neq \0,$ then  we obtain the CCN, denoted by $ \mathfrak{K}^{\infty}_{com,\Upsilon}(\bm v).$
  \end{itemize}
\end{remark}

Next, we present a key result that is essential for the following sections. To derive this, 
let $\mathbf{H}=(A,B,C,D,E)$ and we set $$\vec(\mathbf{H})=[\vec(A)^{\bm\top}, \vec(B)^{\bm\top}, \vec(C)^{\bm\top}, \vec(D)^{\bm\top}, \vec(E)^{\bm\top}]^{\bm\top}.$$
Consider $\D A, \D B, \D C, \D D, \D E,$ and $\D \mathbf{b}=[\D {\bf b}_1^{\T}, \D {\bf b}_2^{\T}, \D {\bf b}_3^{\T}]^{\T}$ are the perturbations on $A, B, C, D,E,$ and $\mathbf{b},$ respectively. Further, we denote
$$\D \mathfrak{B}=\bmatrix{\D A & \D B^{\bm\top} & \bf 0\\ \D B & -\D D & \D C^{\bm\top}\\ \bf 0 & \D C& \D E},$$ and assume that $\|\D \mathfrak{B}\|_2\leq \bm{\epsilon} \|\mathfrak{B}\|_2$ 
and $\|\D {\bf b}\|_2\leq \bm{\epsilon} \|{\bf b}\|_2$. Then, we have the following perturbed DSSP: 
\begin{eqnarray}\label{eq:pert}
    \bmatrix{A+\D A & (B+\D B)^{\bm\top} & \bf 0\\ B+\D B &-(D+ \D D) & (C+ \D C)^{\bm\top} \\ \bf 0 & C+\D C& E +\D E} \bmatrix{\bm{x}+\D \bm{x}\\ \bm{y}+ \D \bm{y} \\ \bm{z} +\D \bm{z}}= \bmatrix{\mathbf{b}_1+\D \mathbf{b}_1\\ \mathbf{b}_2+\D \mathbf{b}_2 \\\mathbf{b}_3+\D \mathbf{b}_3},
\end{eqnarray}
which has the unique solution $\bmatrix{\bm{x}+\D \bm{x}\\ \bm{y}+ \D \bm{y} \\ \bm{z} +\D \bm{z}}$ when $\|\mathfrak{B}^{-1}\|_2\|\Delta \mathfrak{B}\|_2<1.$ 

Consequently, we obtain the following important result.
\begin{lemma}\label{lemma:sec2}
Suppose $[\x^{\T}, \y^{\T}, \z^{\T}]$ and $[(\x+\D \x)^{\T}, (\y+\D \y)^{\T}, (\z+\D \z)^{\T}]^{\T}$ are the unique solutions of the original DSPP \eqref{eq1:SPP} and perturbed DSPP \eqref{eq:pert}, respectively.
Then, the first-order perturbation expression of $\bmatrix{\D\bm{x}^{\T},\D\bm{y}^{\T}, \D\bm{z}^{\T}}^{\T}$ is given by 
\begin{eqnarray}\label{eq27}
  \bmatrix{\D\bm{x}\\ \D\bm{z}\\ \D\bm{y}}= -\mathfrak{B}^{-1} \bmatrix{ \mathcal{G} &  -I_{\bm l}}\bmatrix{\vec(\D \mathbf{H})\\ \D {\bf b}} + \mathcal{O}(\bm{\epsilon}^2),
\end{eqnarray}
 where \begin{equation}
     \mathcal{G}= \bmatrix{\bm{x}^{\bm\top}\otimes I_n & I_n\otimes \bm{y}^{\bm\top}& \bf 0& \bf 0& \bf 0\\ \bf 0& \bm{x}^{\bm\top}\otimes I_m & I_m\otimes \bm{z}^{\bm\top} &- \bm{y}^{\bm\top}\otimes I_m & \bf 0\\ \bf 0 & \bf 0& \bm{y}^{\bm\top}\otimes I_p & \bf 0 & \bm{z}^{\bm\top}\otimes I_p}\in \R^{\bm{l} 
  \times \bm{s}},
 \end{equation}
 $$\vec(\D \mathbf{H})=[\vec(\D A)^{\bm\top}, \vec(\D B)^{\bm\top}, \vec(\D C)^{\bm\top}, \vec(\D D)^{\bm\top}, \vec(\D E)^{\bm\top}]^{\bm\top},$$
 and ${\bm s}=(n^2+m^2+p^2+nm+mp).$
 \end{lemma}
\proof
Combining \eqref{eq1:SPP} and \eqref{eq:pert}, we obtain
\begin{eqnarray}\label{eq29}
    \bmatrix{A & B^{\bm\top} & \bf 0\\ B& -D &C^{\bm\top}\\ \bf 0& C & E} \bmatrix{\D \bm{x}\\ \D \bm{y}\\ \D \bm{z}}= \bmatrix{\D \mathbf{b}_1\\ \D \mathbf{b}_2\\ \D \mathbf{b}_3}-\bmatrix{\D A \bm{x}+\D B^{\bm\top} \bm{y}\\ \D B \bm{x}- \D D \bm{y}+ \D C^{\bm\top}\bm{z}\\ \D C \bm{y}+ \D E \bm{z}}+\mathcal{O}(\bm{\epsilon}^2).
\end{eqnarray}
Thus, the proof follows by applying the vec operator on \eqref{eq29} and utilizing the properties listed in \eqref{kron}.
$\blacksquare$

 \section{ Partial Unified CNs for the DSPP }\label{sec:PCN}
This section primarily focuses on developing a unified framework for the partial CN for the solution $\bm{w}=[\x^{\bm\top},\y^{\bm\top},\z^{\bm\top}]^{\bm\top}$ of the DSPP \eqref{eq1:SPP}. As special cases, we also derive the compact formulae and computationally efficient upper bounds for the partial NCN, MCN and CCN.

To derive the  partial unified CN of the DSPP \eqref{eq1:SPP}, we define the following mapping:
\begin{align}\label{eq:map}
\nonumber \vp:\,& \R^{n\times n}\times  \R^{m\times n}\times  \R^{p\times m}\times  \R^{m\times m}\times  \R^{p\times p}\times \R^{\bm l} \rightarrow \R^k\\
&~~\vp(\mathbf{H}, \mathbf{b})=\mathrm{L}\bmatrix{\x\\ \y \\ \z}=\mathrm{L}\mathfrak{B}^{-1}\mathbf{b},
\end{align}
where $\mathrm{L}\in \R^{k\times \bm l} (k\leq \l).$
Following the Definition \ref{def:ucn}, we now define the partial unified CN for the DSPP using the mapping $\vp$ as follows.


\begin{definition}\label{def:PartialCN}
Suppose $\bm{w}=[\bm{x}^{\bm\top}, \bm{y}^{\bm\top}, \bm{z}^{\bm\top}]^{\bm\top}$ is the unique solution of the DSPP \eqref{eq1:SPP} and $\mathrm{L}\in \R^{k\times \bm{l}}.$ Consider the map $\vp$ defined as in \eqref{eq:map}. Then, the  partial unified  CN of $\bm{w}=[\bm{x}^{\bm\top}, \bm{y}^{\bm\top}, \bm{z}^{\bm\top}]^{\bm\top}$  with respect to (w.r.t.) $\mathrm{L}$ is defined as follows:
    \begin{eqnarray*}
         \mathfrak{K}_{\vp}(\mathbf{H}, \b; \mathrm{L}):=\lim_{\bm{\epsilon}\rightarrow 0} \sup_{0<\left\|\vec\left(\Psi^{\ddagger} \odot \Delta \mathbf{H},\, \be^{\ddagger} \odot \D \b\right)\right\|_{\tau}\leq \bm{\epsilon}} \frac{\left\|\xi^{\ddagger}_{\mathrm{L}}\odot \left( \bm{\varphi}(\mathbf{H}+\D \mathbf{H}, \mathbf{b}+ \D \mathbf{b})-\vp(\mathbf{H}, \b)\right) \right\|_\gamma}{\left\|\vec\left(\Psi^{\ddagger} \odot \Delta \mathbf{H}, \be^{\ddagger} \odot \D \b\right)\right\|_{\tau}},
    \end{eqnarray*}
    where $\xi_{\mathrm{L}}\in \R^{k},$ $\Psi=(\Psi_A, \Psi_B,\Psi_C, \Psi_D,\Psi_E),$ $\Psi_A\in \R^{n\times n},$ $\Psi_B\in \R^{m\times n},$ $\Psi_C\in \R^{p\times m},$ $\Psi_D\in \R^{m\times m},$ $\Psi_E\in \R^{p\times p}$ and $\bm{\chi} \in \R^{\bm l}$ are the parameters with the assumptions that 
    if some entries of $\Psi$ and $\be$ are zero, then the corresponding entry of $\D \mathbf{H}$ and $\D {\bf b},$ respectively,  must be zero.
\end{definition}
\begin{remark}
In the context of  Remark \ref{re21}, to obtain the partial NCN, we consider  
$\xi_{\mathrm{L}}=\|\mathrm{L}[\bm x^{\bm\top}, \bm y^{\bm\top}, \bm z^{\bm\top}]^{\T}\|_2 {\bf 1}_{\bm l},$
for the partial MCN, we consider 
$\xi_{\mathrm{L}}=\|\mathrm{L}[\bm x^{\bm\top}, \bm y^{\bm\top}, \bm z^{\bm\top}]^{\T}\|_{\infty} {\bf 1}_{\bm l},$ and
for partial CCN, we consider $\xi_{\mathrm{L}}=\mathrm{L}[\bm x^{\bm\top}, \bm y^{\bm\top}, \bm z^{\bm\top}]^{\T}.$
\end{remark}
In the following theorem, we provide a compact and closed-form expression for the partial unified CN.
\begin{theorem}\label{th1}
    Suppose $\bm{w}=[\bm{x}^{\bm\top}, \bm{y}^{\bm\top}, \bm{z}^{\bm\top}]^{\bm\top}$ is the unique solution of the DSPP \eqref{eq1:SPP} and $\mathrm{L}\in \R^{k\times \bm{l}}.$ Then, the partial unified CN of $\bm{w}=[\bm{x}^{\bm\top}, \bm{y}^{\bm\top}, \bm{z}^{\bm\top}]^{\bm\top}$ w.r.t. $\mathrm{L}$ is given by
    \begin{equation}
         \mathfrak{K}_{\vp}(\mathbf{H}, {\bf  b}; \mathrm{L}) = \left\|\Theta_{\xi^{\ddagger}_{\mathrm{L}}}\mathrm{L}\mathfrak{B}^{-1}\bmatrix{ \mathcal{G} &  -I_{\bm l}} \bmatrix{\Theta_{\vec(\Psi)}& \bf 0\\ \bf 0 & \Theta_{\be}}\right\|_{\tau,\gamma},
    \end{equation}
    where $\|\cdot~\|_{\tau,\gamma}$ is the matrix norm induced by vector norms $\|\cdot~\|_{\tau} $ and $\|\cdot~\|_{\gamma}.$
\end{theorem}
\proof From the definition of the mapping $\vp$ in \eqref{eq:map} and Lemma \ref{lemma:sec2}, we get
\begin{align}
\nonumber    \vp(\mathbf{H}+\D \mathbf{H}, \b +\D \b)-\vp(\mathbf{H}, \b)&=\mathrm{L}\bmatrix{\bm{x}+\D\bm{x}\\ \bm{y}+\D\bm{y}\\ \bm{z}+\D\bm{z}}-\mathrm{L}\bmatrix{\x\\ \y \\ \z}\\ \nonumber
    &=\mathrm{L} \bmatrix{\D\bm{x}\\ \D\bm{y}\\ \D\bm{z}}\\ \label{th1:eq1}
    &= -\mathrm{L}\mathfrak{B}^{-1}\bmatrix{ \mathcal{G} &  -I_{\bm l}}\bmatrix{\vec(\D \mathbf{H})\\ \D \b}+ \mathcal{O}(\bm{\epsilon}^2).
\end{align}
By considering the requirement on $\Psi$ and $\be,$ we have
\begin{equation}\label{th1:eq2}
    \bmatrix{\vec(\D \mathbf{H}) \\ \D {\bf b}}=\bmatrix{\Theta_{\vec(\Psi)}& \bf 0\\ \bf 0 & \Theta_{\be}}\bmatrix{\vec(\Psi^{\ddagger} \odot \D \mathbf{H}) \\ \bm{\chi}^{\ddagger} \odot \D {\bf b}}.
\end{equation}
Substituting \eqref{th1:eq2} into \eqref{th1:eq1} and from Definition \ref{def:PartialCN}, we obtain
\begin{align}
     \nonumber  \mathfrak{K}_{\vp}(\mathbf{H}, \b; \mathrm{L})& = \sup_{\left\|\vec\left(\Psi^{\ddagger} \odot \Delta \mathbf{H},\, \be^{\ddagger} \odot \D \b\right)\right\|_{\tau}\neq 0} \frac{\left\|\Theta_{\xi^{\ddagger}_{\mathrm{L}}}\mathrm{L} 
 \mathfrak{B}^{-1} \bmatrix{ \mathcal{G} &  -I_{\bm l}} \bmatrix{\Theta_{\vec(\Psi)}& \bf 0\\ \bf 0 & \Theta_{\be}} \bmatrix{\vec(\Psi^{\ddagger} \odot \D \mathbf{H})\\ \bm{\chi}^{\ddagger} \odot \D {\bf b}}\right\|_{\gamma}}{\left\|\vec\left(\Psi^{\ddagger} \odot \Delta \mathbf{H}, \be^{\ddagger} \odot \D \b\right)\right\|_{\tau}}\\
      &=\left\|\Theta_{\xi^{\ddagger}_{\mathrm{L}}}\mathrm{L} \mathfrak{B}^{-1}  \bmatrix{\mathcal{G} &  -I_{\bm l}} \bmatrix{\Theta_{\vec(\Psi)}& \bf 0\\ \bf 0 & \Theta_{\be}}\right\|_{\tau,\gamma}.
\end{align}
Hence, the proof is completed.
$\blacksquare$

Next, we derive various partial CNs by considering specific norms. In the following result, we focus on   when \(\tau = \gamma = 2\).
\begin{theorem}\label{th32}
    Consider $\tau=\gamma= 2$ and assuming that  $\Psi,$ $\be$ and $\xi_{\mathrm{L}}$ are positive real numbers, then the partial CN has the following forms:
    \begin{align}\label{eq1:th2}
          &  \mathfrak{K}^{(2)}_{\vp}(\mathbf{H}, \b; \mathrm{L}) = \frac{\left\|\mathrm{L} \mathfrak{B}^{-1}   \bmatrix{\Psi  \mathcal{G} &  -\be I_{\bm l}} \right\|_{2}}{\xi_{\mathrm{L}}},\\ \label{eq:th32}
         & \widehat{ \mathfrak{K}}^{(2)}_{\vp}(\mathbf{H}, \b; \mathrm{L}) = \frac{\left\|\mathrm{L}\mathfrak{B}^{-1}({\Psi}^2 \mathcal{J}   +{\be}^2I_{\bm l}) (\mathfrak{B}^{-1})^{\T}\mathrm{L}^{\T} \right\|^{1/2}_{2}}{\xi_{\mathrm{L}}},
         %
    \end{align}
    where $\mathcal{J}\in \R^{\bm{l}\times \bm{l}}$ is given by
    $$\mathcal{J}=\bmatrix{(\|\x\|_2^2 +\|\y\|_2^2) I_n & \x \y^{\bm\top} & {\bf 0}\\ \y\x^{\bm\top}  & (\|\x\|_2^2 +\|\y\|_2^2+ \|\z\|_2^2)  I_m & \y\z^{\bm\top}\\ \bf 0  & \z\y^{\bm\top}  & (\|\y\|_2^2 +\|\z\|_2^2)  I_p}.$$
\end{theorem}
\proof Since $\tau=\gamma =2$ and $\Psi,$ $\be$ and $\xi_{\mathrm{L}}$ are positive real numbers, from Theorem \ref{th1}, we obtain \begin{eqnarray}\label{eq3:th2}
     \mathfrak{K}^{(2)}_{1,\vp}(\mathbf{H}, \b; \mathrm{L})= \frac{\left\|\mathrm{L}\mathfrak{B}^{-1}\bmatrix{\Psi \mathcal{G} &  -\be I_{\bm l}} \right\|_{2}}{\xi_{\mathrm{L}}}.
\end{eqnarray}
 Using the property that, for $Z\in \R^{m\times n},$ $\|Z\|_2 =\|ZZ^{\bm\top}\|_2^{1/2},$ we have 
 \begin{align}\label{eq2:th2}
     \nonumber \left\|\mathrm{L}\mathfrak{B}^{-1}\bmatrix{\Psi \mathcal{G} &  -\be I_{\bm l}} \right\|_{2}&= \left\| \mathrm{L}\mathfrak{B}^{-1} (\Psi^2\mathcal{G} \mathcal{G}^{\bm\top}+\be^2 I_{\bm l}) (\mathfrak{B}^{-1})^{\T}\mathrm{L}^{\bm\top}\right\|_2^{1/2}\\
     &= \left\| \mathrm{L}\mathfrak{B}^{-1} (\Psi^2\mathcal{J} +\be^2 I_{\bm l}) (\mathfrak{B}^{-1})^{\T}\mathrm{L}^{\bm\top}\right\|_2^{1/2}.
 \end{align}
 Hence, the proof follows by substituting \eqref{eq2:th2} into  \eqref{eq3:th2}.
$\blacksquare$
\begin{remark}
  The partial CN in Theorem \ref{th32} is a simplified version of the partial NCN for the solution of the DSPP \eqref{eq1:SPP}. The  NCN for ${\bm w}=[\x^{\T},  \y^{\T},\z^{\T  }]^{\T}, \x, \y$ and $\z$ can be obtained by considering $$\mathrm{L}= I_{\bm l},~ \bmatrix{I_n & {\bf 0}_{n\times (m+p)}}, ~\bmatrix{{\bf 0}_{m\times n}& I_m&{\bf 0}_{m\times p}} ~\text{and}~ \bmatrix{{\bf 0}_{p\times (n+m)}& I_p},$$ respectively, in Theorem \ref{th32}.
\end{remark}
\begin{remark}\label{re3}
  {Theorem \ref{th32} presents two equivalent expressions \eqref{eq1:th2} and \eqref{eq:th32}  of the partial NCN. Note that the sizes of the matrices present in these expressions are   \(k \times (\bm{l}+n^2+m^2+p^2+nm+mp)\) and \(k \times k\), respectively. Hence, the expression in \eqref{eq:th32} significantly reduces the storage requirements. Moreover, equation \eqref{eq1:th2} contains Kronecker products, which significantly increase memory requirements and computational complexity, making it less practical for large-scale problems. In contrast, equation \eqref{eq:th32} avoids Kronecker products and is more efficient, especially when $k< \l$, since it results in smaller matrix sizes. 
  Further, to improve computational efficiency and avoid computing $\mathfrak{B}^{-1}$ explicitly present in the expressions of Theorem \ref{th32}, we adopt a practical approach inspired by the procedure outlined in  \cite{UCN2018}. For \eqref{eq1:th2}, we first solve $\mathfrak{B}X=\bmatrix{\Psi  \mathcal{G} &  -\be I_{\bm l}}$  by using LU decomposition, and then compute $\mathrm{L}X.$ For the expression in \eqref{eq:th32}, we first solve the system $\mathfrak{B}Y={\Psi}^2 \mathcal{J}   +{\be}^2I_{\bm l},$ and following that solve the system $\mathfrak{B}^{\T}Z=Y,$ and compute $\mathrm{L} Z \mathrm{L}^{\T}.$ Note that we need to perform only one LU decomposition of $\mathfrak{B}.$ 
 }
\end{remark}

In the next result,  {to address the computational challenges mentioned in Remark \ref{re3} associated with the expressions in Theorem~\ref{th32},} we provide an easily computable upper bound for the partial NCN $   \mathfrak{K}^{(2)}_{\Phi}(\mathbf{H}, \b; \mathrm{L}).$
\begin{corollary}\label{cor31}
    Under the assumption of Theorem \ref{th32}, we have the following upper bound:
    \begin{eqnarray}
         \mathfrak{K}^{(2)}_{1,\vp}(\mathbf{H}, \b; \mathrm{L})\leq \mathfrak{K}^{(2),u}_{1,\vp}(\mathbf{H}, \b; \mathrm{L}):=\frac{\|\mathrm{L}\mathfrak{B}^{-1}\|_2}{\xi_{\mathrm{L}}}\left( \Psi\|\mathcal{J}\|^{1/2}_{2}+\be \right).
    \end{eqnarray}
\end{corollary}
\proof Using the properties of the spectral norm that for the matrices $X$ and $Y$ of appropriate sizes, $\left\|\bmatrix{X & Y} \right\|_2\leq \|X\|_2 + \|Y\|_2$ and $\|XY\|_2\leq \|X\|_2 \|Y\|_2,$ and from \eqref{eq1:th2}, we obtain 
\begin{align}\label{eq:coro31}
  \nonumber  \left\|\mathrm{L}\mathfrak{B}^{-1}\bmatrix{\Psi \mathcal{G} &  -\be I_{\bm l}} \right\|_{2} &\leq \Psi \| \mathrm{L}\mathfrak{B}^{-1}\mathcal{G}\|_2 +\be \| \mathrm{L}\mathfrak{B}^{-1}\|_2\\ \nonumber
    & \leq \Psi \| \mathrm{L}\mathfrak{B}^{-1}\|_2 \|\mathcal{G}\|_2 +\be \| \mathrm{L}\mathfrak{B}^{-1}\|_2\\
    &= \Psi \| \mathrm{L}\mathfrak{B}^{-1}\|_2 \|\mathcal{J}\|^{1/2}_2 +\be \| \mathrm{L}\mathfrak{B}^{-1}\|_2.
\end{align}
Hence, the proof follows from \eqref{eq1:th2} and \eqref{eq:coro31}. $\blacksquare$

In the following theorem, we investigate the partial CN for the DSPP when \(\tau = \gamma = \infty\), from which we derive the partial MCN and CCN.
\begin{theorem}\label{th33}
    When $\tau=\gamma=\infty,$ the partial CN is given as follows:
    \begin{align}
         \mathfrak{K}^{\infty}_{\vp}(\mathbf{H}, \b; \mathrm{L}) = \left\| |\Theta_{\xi^{\ddagger}_{\mathrm{L}}}| \left|\mathrm{L}\mathfrak{B}^{-1} \bmatrix{ \mathcal{G} &  -I_{\bm l}}\right| \bmatrix{\vec(|\Psi|)\\  |\be|}\right\|_{\infty}.
    \end{align}
    Moreover, we consider $\Psi=\mathbf{H},$ $\be ={\bf b}.$ Then, if we set $\xi_{\mathrm{L}}=\|\mathrm{L}[\bm{x}^{\bm\top}, \bm{y}^{\bm\top}, \bm{z}^{\bm\top}]^{\bm\top}\|_{\infty}{\bf{1}}_{\bm l},$  the partial MCN 
  is given by 
   \begin{align}\label{PMCN}
         \mathfrak{K}^{\infty}_{mix,  \vp}(\mathbf{H}, \b; \mathrm{L}) = \frac{\left\| \left|\mathrm{L}\mathfrak{B}^{-1}\bmatrix{ \mathcal{G} &  -I_{\bm l}}\right| \bmatrix{\vec(|\mathbf{H}|)\\ |{\bf b}| }\right\|_{\infty}}{\|\mathrm{L}[\bm{x}^{\bm\top}, \bm{y}^{\bm\top}, \bm{z}^{\bm\top}]^{\bm\top}\|_{\infty}}
    \end{align}
    and if we set $\xi_{\mathrm{L}}= \mathrm{L}[\bm{x}^{\bm\top}, \bm{y}^{\bm\top}, \bm{z}^{\bm\top}]^{\bm\top},$  the partial CCN is given by
    \begin{align}\label{PCCN}
         \mathfrak{K}^{\infty}_{com, \vp}(\mathbf{H}, \b; \mathrm{L}) = \left\|\frac{ \left|\mathrm{L} \mathfrak{B}^{-1} \bmatrix{ \mathcal{G} &  -I_{\bm l}}\right| \bmatrix{\vec(|\mathbf{H}|)\\ |{\bf b}| }}{\left|\mathrm{L}[\bm{x}^{\bm\top}, \bm{y}^{\bm\top}, \bm{z}^{\bm\top}]^{\bm\top}\right|}\right\|_{\infty}.
    \end{align}
\end{theorem}
\proof Consider $\tau=\gamma=\infty,$ then from Theorem \ref{th1}, we have
\begin{align*}
     \mathfrak{K}^{\infty}_{\Phi}(\mathbf{H}, \b; \mathrm{L})& = \left\| \Theta_{\xi^{\ddagger}_{\mathrm{L}}}\mathrm{L}\mathfrak{B}^{-1} \bmatrix{ \mathcal{G} &  -I_{\bm l}} \bmatrix{\Theta_{\vec(\Psi)}& \bf 0\\ \bf 0 & \Theta_{\be}}\right\|_{\infty}\\
    &= \left\| |\Theta_{\xi^{\ddagger}_{\mathrm{L}}}|  \left|\mathfrak{B}^{-1} \bmatrix{ \mathcal{G} &  -I_{\bm l}}\right| \bmatrix{|\Theta_{\vec(\Psi)}|& \bf 0\\ \bf 0 & |\Theta_{\be}|} \right\|_{\infty}\\
     &= \left\| |\Theta_{\xi^{\ddagger}_{\mathrm{L}}}|  \left|\mathfrak{B}^{-1} \bmatrix{ \mathcal{G} &  -I_{\bm l}}\right| \bmatrix{|\Theta_{\vec(\Psi)}|& \bf 0\\ \bf 0 & |\Theta_{\be}|} {\bf 1}_{s+ \bm{l}}\right\|_{\infty}\\
     &= \left\| |\Theta_{\xi^{\ddagger}_{\mathrm{L}}}|  \left| \mathfrak{B}^{-1} \bmatrix{ \mathcal{G} &  -I_{\bm l}}\right| \bmatrix{\vec(|\Psi|)\\  |\be|} \right\|_{\infty}.
 \end{align*}
Rest of the proof followings considering $\Psi = \mathbf{H},$ $\be= {\bf b},$ $\xi_{\mathrm{L}}=\|\mathrm{L}[\bm{x}^{\bm\top}, \bm{y}^{\bm\top}, \bm{z}^{\bm\top}]^{\bm\top}\|_{\infty}{\bf{1}}_{\bm l} $ or $\xi_{\mathrm{L}}=\mathrm{L}[\bm{x}^{\bm\top}, \bm{y}^{\bm\top}, \bm{z}^{\bm\top}]^{\bm\top}.$
$\blacksquare$

\begin{remark}
  The MCN and CCN for ${\bm w}=[\x^{\T},  \y^{\T},\z^{\T  }]^{\T}, \x, \y$ and $\z$ can be obtained by considering $$\mathrm{L}= I_{\bm l},~ \bmatrix{I_n & {\bf 0}_{n\times (m+p)}}, ~\bmatrix{{\bf 0}_{m\times n}& I_m&{\bf 0}_{m\times p}} ~\text{and}~ \bmatrix{{\bf 0}_{p\times (n+m)}& I_p},$$ respectively, in \eqref{PMCN}  and \eqref{PCCN} of Theorem \ref{th33}.
\end{remark}
In the following result, we provide sharp, computationally less expensive and Kronecker product free upper bounds for the partial MCN and CCN obtained in Theorem \ref{th33}.
\begin{corollary}\label{coro32}
    Assume that the conditions in Theorem \ref{th33} hold. Then
    \begin{align*}
         \mathfrak{K}^{\infty}_{mix, \vp}(\mathbf{H}, \b; \mathrm{L})\leq  \mathfrak{K}^{\infty, u}_{mix, \vp}(\mathbf{H}, \b; \mathrm{L}) := \frac{\left\||\mathrm{L}\mathfrak{B}^{-1}| \left(|\mathcal{H}|+|{\bf b}|\right)\right\|_{\infty}}{\left\|\mathrm{L}[\bm{x}^{\bm\top}, \bm{y}^{\bm\top}, \bm{z}^{\bm\top}]^{\bm\top}\right\|_{\infty}}
    \end{align*}
    and \begin{align*}
         \mathfrak{K}^{\infty}_{com, \vp}(\mathbf{H}, \b; \mathrm{L})\leq  \mathfrak{K}^{\infty, u}_{com, \vp}(\mathbf{H}, \b; \mathrm{L}) : = \left\|\frac{|\mathrm{L}\mathfrak{B}^{-1}| \left(|\mathcal{H}|+|{\bf b}|\right)}{\mathrm{L}[\bm{x}^{\bm\top}, \bm{y}^{\bm\top}, \bm{z}^{\bm\top}]^{\bm\top}}\right\|_{\infty},
    \end{align*}
    where $\mathcal{H}=\bmatrix{|A||\x|+|B^{\bm\top}||\y|\\ |B||\x|+|D^{\bm\top}||\y|+|C^{\bm\top}||\z|\\ |C||\y|+|E||\z|}.$
\end{corollary}
\proof Utilizing the properties of Kronecker product in \eqref{kron}, we have 
\begin{align}\label{eq:coro1}
  \nonumber  &\left|\mathrm{L}\mathfrak{B}^{-1}\bmatrix{ \mathcal{G} ~&  -I_{\bm l}}\right| \bmatrix{\vec(|\mathbf{H}|)\\ |{\bf b}| }\leq |\mathrm{L} \mathfrak{B}^{-1}| \bmatrix{|\mathcal{G}| & I_{\bm l}} \bmatrix{\vec(|\mathbf{H}|)\\ |{\bf b}| } \\ \nonumber
   &= |\mathrm{L} \mathfrak{B}^{-1}| (|\mathcal{G}| \vec(|\mathbf{H}|) + |{\bf b}| )\\ \nonumber
   &=  |\mathrm{L} \mathfrak{B}^{-1}| \left(\bmatrix{ (|\x|^{\bm\top} \otimes I_n) \vec(|A|) + (I_n\otimes |\y|^{\bm\top}) \vec(|B|)\\ (|\x|^{\bm\top} \otimes I_m) \vec(|B|)+ (I_m\otimes |\z|^{\bm\top}) \vec(|C|) + (|\y|^{\bm\top} \otimes I_m) \vec(|D|) \\ (|\y|^{\bm\top} \otimes I_p) \vec(|C|) + (|\z|^{\bm\top} \otimes I_p) \vec(|E|)} + |{\bf b}|\right)\\  
   &= |\mathrm{L}\mathfrak{B}^{-1}| \left(\left|\bmatrix{|A||\x|+|B^{\bm\top}||\y|\\ |B||\x|+|D^{\bm\top}||\y|+|C^{\bm\top}||\z|\\ |C||\y|+|E||\z|}\right| + |{\bf b}|\right).
\end{align}
From \eqref{eq:coro1} and the expressions of partial MCN and CCN in Theorem \ref{th33}, we get 
\begin{eqnarray*}
    \mathfrak{K}^{\infty}_{mix, \vp}(\mathbf{H}, \b; \mathrm{L})\leq  \mathfrak{K}^{\infty, u}_{mix, \vp}(\mathbf{H}, \b; \mathrm{L}) ~ \text{and}~    \mathfrak{K}^{\infty}_{com, \vp}(\mathbf{H}, \b; \mathrm{L})\leq  \mathfrak{K}^{\infty, u}_{com, \vp}(\mathbf{H}, \b; \mathrm{L}).
\end{eqnarray*}
Hence, the proof follows.
$\blacksquare$

\section{Structured partial CNs}\label{sec:SPCN}
Consider three subspaces $\mathbb{S}_1\subseteq \R^{n\times n}$, $\mathbb{S}_2\subseteq \R^{m\times m}$ and 
     $\mathbb{S}_3\subseteq \R^{p\times p}$ consisting of three distinct linear structured matrices, such as symmetric. Suppose that the corresponding dimensions of the linear subspaces are $s,$ $r$, and $q,$ respectively. Let $A\in \mathbb{S}_1, $ $D\in \mathbb{S}_2$ and $E\in \mathbb{S}_3,$ then according to \cite{HIGHAM1992, TLS2011Li, Rump2003}, there exist unique generating vectors $\bm{a}\in \R^{s}$ $\bm{d}\in \R^{r}$ and $\bm{e}\in \R^{q}$ such that 
\begin{align}\label{eq:51}
    \vec(A)=\bm{\Phi}_{\mathbb{S}_1}\bm{a}, ~\vec(D)=\bm{\Phi}_{\mathbb{S}_2}\bm{d} ~~\text{and}~~ \vec(E)=\bm{\Phi}_{\mathbb{S}_3}\bm{e},
\end{align}
where $\bm{\Phi}_{\mathbb{S}_1}\in \R^{n^2\times s},$ $\bm{\Phi}_{\mathbb{S}_2}\in \R^{m^2\times r}$ and $\bm{\Phi}_{\mathbb{S}_3}\in \R^{p^2\times q}.$ These matrices are fixed for each specific structure and encapsulate the information corresponding to the linear structure of their respective subspaces.  

Let $\vec_{\mathbb{S}}(\mathbf{H})=[\bm{a}^{\bm\top}, \vec(B)^{\bm\top}, \vec(C)^{\bm\top}, \bm{d}^{\bm\top}, \bm{e}^{\bm\top}]^{\bm\top}.$ Then the  structured partial CN for solution $\bm{w}=[\x^{\bm\top},\y^{\bm\top},\z^{\bm\top}]^{\bm\top}$ of the DSPP \eqref{eq1:SPP}  w.r.t. $\mathrm{L}$ is given by
\begin{align}\label{def:SCN}
         \mathfrak{K}^{\mathbb{S}}_{\vp}(\mathbf{H}, \b; \mathrm{L}):=\lim_{\bm{\epsilon}\rightarrow 0} \sup_{\substack{0<\left\|\vec\left(\Psi^{\ddagger} \odot \Delta \mathbf{H},\, \be^{\ddagger} \odot \D \b\right)\right\|_{\tau}\leq \bm{\epsilon}\\ \D A\in\, \mathbb{S}_1,\, \D D\in\, \mathbb{S}_2,\, \D E \in\, \mathbb{S}_3} } \frac{\left\|\xi^{\ddagger}_{\mathrm{L}}\odot \left( \vp(\mathbf{H}+\D \mathbf{H}, \mathbf{b}+ \D \mathbf{b})-\vp(\mathbf{H}, \b)\right) \right\|_\gamma}{\left\|\vec\left(\Psi^{\ddagger} \odot \Delta \mathbf{H}, \be^{\ddagger} \odot \D \b\right)\right\|_{\tau}},
    \end{align}
    where $\xi_{\mathrm{L}}\in \R^{k},$ $\Psi=(\Psi_A, \Psi_B,\Psi_C, \Psi_D,\Psi_E),$ $\Psi_A\in \mathbb{S}_1,$ $\Psi_B\in \R^{m\times n},$ $\Psi_C\in \R^{p\times m},$ $\Psi_D\in \mathbb{S}_2,$ $\Psi_E\in \mathbb{S}_3$ and $\bm{\chi} \in \R^{\bm l}$.

    Since the matrices $\D A, \Psi_A\in \mathbb{S}_1, $ $\D D, \Psi_D\in \mathbb{S}_2$ and $\D E, \Psi_E\in \mathbb{S}_3,$ as in \eqref{eq:51}, we have 
    \begin{align}\label{eq:53}
  &  \vec(\D A)=\bm{\Phi}_{\mathbb{S}_1}\D \bm{a}, ~ \vec(\Psi_A)=\bm{\Phi}_{\mathbb{S}_1}\psi_A, ~\vec(\D D)=\bm{\Phi}_{\mathbb{S}_2}\D\bm{d}, \\
    & \vec(\Psi_D)=\bm{\Phi}_{\mathbb{S}_2}\psi_D,~ \vec(\D E)=\bm{\Phi}_{\mathbb{S}_3}\D \bm{e}, ~ \vec(\Psi_E)=\bm{\Phi}_{\mathbb{S}_3}\psi_E,
\end{align}
where $\D \bm{a},$  $\D \bm{d},$  $\D \bm{e},$ $\psi_D,$ $\psi_A,$ and $\psi_E$ are the unique generating vectors of $\D A,$ $\D D,$ $\D E,$ $\Psi_A,$  $\Psi_D,$  and $\Psi_E$, respectively. Note that, $\Psi_A^{\ddagger}\in \mathbb{S}_1,$ $\Psi_D^{\ddagger}\in \mathbb{S}_2$ and $\Psi_E^{\ddagger}\in \mathbb{S}_3$. Consequently, we obtain
\begin{align}\label{eq:55}
    \vec(\Psi_A^{\ddagger}\odot\D A)&=\bm{\Phi}_{\mathbb{S}_1} (\psi_A^{\ddagger}\odot \D\bm{a}),~~ \vec(\Psi_D^{\ddagger}\odot\D D)=\bm{\Phi}_{\mathbb{S}_2} (\psi_D^{\ddagger}\odot \D\bm{d})\\ \label{eq:56}
     &\vec(\Psi_E^{\ddagger}\odot\D E)=\bm{\Phi}_{\mathbb{S}_3} (\psi_E^{\ddagger}\odot \D\bm{e}).
\end{align}
Subsequently, we obtain the following result.
\begin{lemma}
Let $\D A, \Psi_A \in \mathbb{S}_1,$ $\D D, \Psi_D \in \mathbb{S}_2,$ $\D E, \Psi_E \in \mathbb{S}_3,$ $B, \Psi_B\in \R^{m\times n},$ $C, \Psi_C\in \R^{p\times m},$  and ${\bf b}, \be\in \R^{\bm l}.$ Then, we have 
    \begin{equation}\label{eq:57}
        \bmatrix{\vec(\Psi^{\ddagger} \odot \mathbf{H})\\ \be^{\ddagger} \odot \D {\bf b}}=\bmatrix{\bm{\Phi}_{\mathbb{S}} & \bf 0 \\ \bf 0 & I_{\bm l}}\bmatrix{\vec_{\mathbb{S}}(\Psi^{\ddagger} \odot H)\\ \be^{\ddagger} \odot \D {\bf b}},
    \end{equation}
    where $\vec_{\mathbb{S}}(\Psi^{\ddagger} \odot \mathbf{H})=[(\psi_A^{\ddagger} \odot \D\bm{a})^{\bm\top}, (\Psi_B^{\ddagger}\odot \vec(\D B))^{\bm\top}, (\Psi_C^{\ddagger}\odot\vec(\D C))^{\bm\top} , (\psi_D^{\ddagger}\odot \D\bm{d})^{\bm\top} , (\psi_D^{\ddagger}\odot  \D\bm{e})^{\bm\top}]^{\bm\top}$ and
    
    $$\bm{\Phi}_{\mathbb{S}}=\bmatrix{\bm{\Phi}_{\mathbb{S}_1} & \bf 0 & \bf 0 & \bf 0\\\bf 0 & I_{mn+mp} & \bf 0 & \bf 0\\ \bf 0  & \bf 0 &\bm{\Phi}_{\mathbb{S}_2} & \bf 0 \\ \bf 0  & \bf 0  & \bf 0 &\bm{\Phi}_{\mathbb{S}_3} }.$$
\end{lemma}
\proof The proof follows immediately using identities in \eqref{eq:55} and \eqref{eq:56}. $\blacksquare$

In the following theorem, we present closed-form expressions for the structured partial CN by considering $\tau=\gamma=2.$
\begin{theorem}\label{Th:SNCN}
    Let $A\in \mathbb{S}_1,$ $D\in \mathbb{S}_2,$ $E\in \mathbb{S}_3$ and $\mathrm{L}\in \R^{k\times \bm{l}}.$ Suppose that ${\bm w}=[\x^{\bm\top},\y^{\bm\top},\z^{\bm\top}]^{\bm\top}$ is the unique solution of DSPP \eqref{eq1:SPP} and consider $\tau=\gamma=2.$ Then the structured partial CN of ${\bm w}=[\x^{\bm\top},\y^{\bm\top},\z^{\bm\top}]^{\bm\top}$ w.r.t. $\mathrm{L}$ is given by \begin{eqnarray*}
\mathfrak{K}^{(2),\mathbb{S}}_{\vp}(\mathbf{H}, \b; \mathrm{L})= \left\|\Theta_{\xi^{\ddagger}_{\mathrm{L}}}\mathrm{L}\mathfrak{B}^{-1}\bmatrix{ \mathcal{G} &  -I_{\bm l}} \bmatrix{\Theta_{\vec(\Psi)}\bm{\Phi}_{\mathbb{S}}\mathfrak{D}^{-1}_{\mathbb{S}}& \bf 0\\ \bf 0 & \Theta_{\be}}\right\|_2,
    \end{eqnarray*}
    where $$\mathfrak{D}_{\mathbb{S}}=\bmatrix{\Theta_{\bm{u}_1} & \bf 0 & \bf 0 & \bf 0\\ \bf 0 & I_{mn+mp}& \bf 0 & \bf 0 \\  \bf 0& \bf 0 & \Theta_{\bm{u}_2}& \bf 0 \\ \bf 0& \bf 0& \bf 0& \Theta_{\bm{u}_3} },$$
    $$\bm{u}_1=[\|\bm{\Phi}_{\mathbb{S}_1}(:,1)\|_2,\ldots, \|\bm{\Phi}_{\mathbb{S}_1}(:,s)\|_2 ]^{\bm\top},~~\bm{u}_2=[\|\bm{\Phi}_{\mathbb{S}_2}(:,1)\|_2,\ldots, \|\bm{\Phi}_{\mathbb{S}_2}(:,r)\|_2 ]^{\bm\top},$$
$$\text{and}~~\bm{u}_3=[\|\bm{\Phi}_{\mathbb{S}_3}(:,1)\|_2,\ldots, \|\bm{\Phi}_{\mathbb{S}_3}(:,q)\|_2 ]^{\bm\top}.$$
\end{theorem}
\proof Taking $\tau=\gamma=2$ in \eqref{def:SCN}, and using \eqref{th1:eq2} and \eqref{th1:eq1}, we obtain
\begin{eqnarray*}
         \mathfrak{K}^{(2),\mathbb{S}}_{\vp}(\mathbf{H}, \b; \mathrm{L})= \sup_{\substack{\left\|\vec\left(\Psi^{\ddagger} \odot \Delta \mathbf{H},\, \be^{\ddagger} \odot \D \b\right)\right\|_{2}\neq 0\\ \D A\in\, \mathbb{S}_1,\, \D D\in\, \mathbb{S}_2,\, \D E \in\, \mathbb{S}_3} } \frac{\left\|\Theta_{\xi^{\ddagger}_{\mathrm{L}}}\mathrm{L} 
 \mathfrak{B}^{-1} \bmatrix{ \mathcal{G} &  -I_{\bm l}} \bmatrix{\Theta_{\vec(\Psi)}& \bf 0\\ \bf 0 & \Theta_{\be}} \bmatrix{\vec(\Psi^{\ddagger} \odot \D \mathbf{H})\\ \bm{\chi}^{\ddagger} \odot \D {\bf b}}\right\|_2}{\left\|\vec\left(\Psi^{\ddagger} \odot \Delta \mathbf{H}, \be^{\ddagger} \odot \D \b\right)\right\|_{2}}.
    \end{eqnarray*}
    Substituting \eqref{eq:57} into the above equation yields
  {\footnotesize   \begin{equation}\label{eq:58}
         \mathfrak{K}^{(2),\mathbb{S}}_{\vp}(\mathbf{H}, \b; \mathrm{L})= \sup_{\substack{\footnotesize\left\|\bmatrix{\bm{\Phi}_{\mathbb{S}} & \bf 0 \\ \bf 0 & I_{\bm l}}\bmatrix{\vec_{\mathbb{S}}(\Psi^{\ddagger} \odot \D\mathbf{H})\\ \be^{\ddagger} \odot \D {\bf b}}\right\|_{2}\neq 0\\ \D A\in\, \mathbb{S}_1,\, \D D\in\, \mathbb{S}_2,\, \D E \in\, \mathbb{S}_3} } \frac{\left\|\Theta_{\xi^{\ddagger}_{\mathrm{L}}}\mathrm{L} 
 \mathfrak{B}^{-1} \bmatrix{ \mathcal{G} &  -I_{\bm l}} \bmatrix{\Theta_{\vec(\Psi)}\bm{\Phi}_{\mathbb{S}} & \bf 0\\ \bf 0 & \Theta_{\be}} \bmatrix{\vec_{\mathbb{S}}(\Psi^{\ddagger} \odot \D\mathbf{H})\\ \be^{\ddagger} \odot \D {\bf b}}\right\|_2}{\left\|\bmatrix{\bm{\Phi}_{\mathbb{S}} & \bf 0 \\ \bf 0 & I_{\bm l}}\bmatrix{\vec_{\mathbb{S}}(\Psi^{\ddagger} \odot \D\mathbf{H})\\ \be^{\ddagger} \odot \D {\bf b}}\right\|_{2}}.
    \end{equation}
    }
    
\noindent    Utilizing the fact that the matrices $\bm{\Phi}_{\mathbb{S}_i}$ for $i=1,2,3,$   are column orthogonal \cite{TLS2011Li}, we get $\bm{\Phi}_{\mathbb{S}_i}^{\bm\top}\bm{\Phi}_{\mathbb{S}_i}=\Theta_{\bm{u}_i}^2,$
where $\Theta_{\bm{u}_i}$ for $i=1,2,3,$ are the diagonal matrices.  
Then
\begin{align}\label{eq:49}
   \nonumber \left\|\bmatrix{\bm{\Phi}_{\mathbb{S}} & \bf 0 \\ \bf 0 & I_{\bm l}}\bmatrix{\vec_{\mathbb{S}}(\Psi^{\ddagger} \odot \D\mathbf{H})\\ \be^{\ddagger} \odot \D {\bf b}}\right\|_{2}&=\left\|\bmatrix{\vec_{\mathbb{S}}(\Psi^{\ddagger} \odot \D\mathbf{H})\\ \be^{\ddagger} \odot \D {\bf b}}^{\bm\top}\bmatrix{\bm{\Phi}_{\mathbb{S}}^{\bm\top}\bm{\Phi}_{\mathbb{S}} & \bf 0 \\ \bf 0 & I_{\bm l}}\bmatrix{\vec_{\mathbb{S}}(\Psi^{\ddagger} \odot \D\mathbf{H})\\ \be^{\ddagger} \odot \D {\bf b}}\right\|_2^{1/2}\\
    &= \left\|\bmatrix{\mathfrak{D}_{\mathbb{S}} & \bf 0 \\ \bf 0 & I_{\bm l}}\bmatrix{\vec_{\mathbb{S}}(\Psi^{\ddagger} \odot \D\mathbf{H})\\ \be^{\ddagger} \odot \D {\bf b}}\right\|_2. 
\end{align}
    Observe that 
    \begin{align}\label{eq:59}
\bmatrix{\bm{\Phi}_{\mathbb{S}} & \bf 0 \\ \bf 0 & I_{\bm l}}\bmatrix{\vec_{\mathbb{S}}(\Psi^{\ddagger} \odot \D\mathbf{H})\\ \be^{\ddagger} \odot \D {\bf b}}= \bmatrix{\bm{\Phi}_{\mathbb{S}}\mathfrak{D}^{-1}_{\mathbb{S}} & \bf 0 \\ \bf 0 & I_{\bm l}}\bmatrix{\mathfrak{D}_{\mathbb{S}} & \bf 0 \\ \bf 0 & I_{\bm l}} \bmatrix{\vec_{\mathbb{S}}(\Psi^{\ddagger} \odot \D \mathbf{H})\\ \be^{\ddagger} \odot \D {\bf b}}.
    \end{align}
    Therefore, substituting \eqref{eq:49} and \eqref{eq:59} in \eqref{eq:58}, we obtain
    {\footnotesize \begin{align*}\label{eq:510}
         \mathfrak{K}^{(2),\mathbb{S}}_{\vp}(\mathbf{H}, \b; \mathrm{L})&= \sup_{\substack{\footnotesize\left\|\bmatrix{\mathfrak{D}_{\mathbb{S}}\vec_{\mathbb{S}}(\Psi^{\ddagger} \odot \D\mathbf{H})\\ \be^{\ddagger} \odot \D {\bf b}}\right\|_{2}\neq 0\\ \D A\in\, \mathbb{S}_1,\, \D D\in\, \mathbb{S}_2,\, \D E \in\, \mathbb{S}_3} } \frac{\left\|\Theta_{\xi^{\ddagger}_{\mathrm{L}}}\mathrm{L} 
 \mathfrak{B}^{-1} \bmatrix{ \mathcal{G} &  -I_{\bm l}} \bmatrix{\Theta_{\vec(\Psi)}\bm{\Phi}_{\mathbb{S}}\mathfrak{D}_{\mathbb{S}}^{-1} & \bf 0\\ \bf 0 & \Theta_{\be}} \bmatrix{\mathfrak{D}_{\mathbb{S}}\vec_{\mathbb{S}}(\Psi^{\ddagger} \odot \D \mathbf{H})\\ \be^{\ddagger} \odot \D {\bf b}}\right\|_2}{\left\|\bmatrix{\mathfrak{D}_{\mathbb{S}}\vec_{\mathbb{S}}(\Psi^{\ddagger} \odot \D \mathbf{H})\\ \be^{\ddagger} \odot \D {\bf b}}\right\|_{2}}\\
 &=\normalsize\left\|\Theta_{\xi^{\ddagger}_{\mathrm{L}}}\mathrm{L} 
 \mathfrak{B}^{-1} \bmatrix{ \mathcal{G} &  -I_{\bm l}} \bmatrix{\Theta_{\vec(\Psi)}\bm{\Phi}_{\mathbb{S}}\mathfrak{D}_{\mathbb{S}}^{-1} & \bf 0\\ \bf 0 & \Theta_{\be}} \right\|_2.
    \end{align*}
    }
    Hence, the proof is completed.
$\blacksquare$

Next, we consider $\tau=\gamma=\infty,$ and derive the structured partial MCN and CCN for the DSPP.
\begin{theorem}\label{Th:MCN}
    Let $A\in \mathbb{S}_1,$ $D\in \mathbb{S}_2,$ $E\in \mathbb{S}_3$ and $\mathrm{L}\in \R^{k\times \bm{l}}.$ Suppose that ${\bm w}=[\x^{\bm\top},\y^{\bm\top},\z^{\bm\top}]^{\bm\top}$ is the unique solution of DSPP \eqref{eq1:SPP}.   When $\tau=\gamma=\infty,$ the structured partial CN of the solution ${\bm w}=[\x^{\bm\top},\y^{\bm\top},\z^{\bm\top}]^{\bm\top}$ 
 w.r.t. $\mathrm{L}$ is given as follows:
    \begin{align}
         \mathfrak{K}^{\infty,\,\mathbb{S}}_{\vp}(\mathbf{H}, \b; \mathrm{L}) = \left\| |\Theta_{\xi^{\ddagger}_{\mathrm{L}}}| \left|\mathrm{L}\mathfrak{B}^{-1} \bmatrix{ \mathcal{G} &  -I_{\bm l}} \bmatrix{\bm{\Phi}_{\mathbb{S}} & \bf 0\\ \bf 0& I_{\bm l}}\right| \bmatrix{\vec_{\mathbb{S}}(|\Psi|)\\  |\be|}\right\|_{\infty}.
    \end{align}
    Moreover, we consider $\Psi=\mathbf{H},$ $\be ={\bf b}.$ Then, if we set $\xi_{\mathrm{L}}=\|\mathrm{L}[\bm{x}^{\bm\top}, \bm{y}^{\bm\top}, \bm{z}^{\bm\top}]^{\bm\top}\|_{\infty}{\bf{1}}_{\bm l},$  the structured  partial MCN 
  is given by 
   \begin{align}
         \mathfrak{K}^{\infty,\,\mathbb{S}}_{mix,  \vp}(\mathbf{H}, \b; \mathrm{L}) = \frac{\left\| \left|\mathrm{L}\mathfrak{B}^{-1}\bmatrix{ \mathcal{G} &  -I_{\bm l}} \bmatrix{\bm{\Phi}_{\mathbb{S}} & \bf 0\\ \bf 0& I_{\bm l}} \right| \bmatrix{\vec_{\mathbb{S}}(|\mathbf{H}|)\\ |{\bf b}| }\right\|_{\infty}}{\|\mathrm{L}[\bm{x}^{\bm\top}, \bm{y}^{\bm\top}, \bm{z}^{\bm\top}]^{\bm\top}\|_{\infty}}
    \end{align}
    and if we set $\xi_{\mathrm{L}}= \mathrm{L}[\bm{x}^{\bm\top}, \bm{y}^{\bm\top}, \bm{z}^{\bm\top}]^{\bm\top},$  the  structured partial CCN is given by
    \begin{align}
         \mathfrak{K}^{\infty,\, \mathbb{S}}_{com, \vp}(\mathbf{H}, \b; \mathrm{L}) = \left\|\frac{ \left|\mathrm{L} \mathfrak{B}^{-1} \bmatrix{ \mathcal{G} &  -I_{\bm l}} \bmatrix{\bm{\Phi}_{\mathbb{S}} & \bf 0\\ \bf 0& I_{\bm l}} \right| \bmatrix{\vec_{\mathbb{S}}(|\mathbf{H}|)\\ |{\bf b}| }}{\left|\mathrm{L}[\bm{x}^{\bm\top}, \bm{y}^{\bm\top}, \bm{z}^{\bm\top}]^{\bm\top}\right|}\right\|_{\infty}.
    \end{align}
\end{theorem}
\proof By construction of the matrices $\bm{\Phi}_{\mathbb{S}_1},$ $\bm{\Phi}_{\mathbb{S}_2}$ and $\bm{\Phi}_{\mathbb{S}_3},$ they have almost one nonzero element in each row. Thus, we get
\begin{equation}\label{eq:513}
 \left\| \bmatrix{\vec(\Psi^{\ddagger} \odot\D \mathbf{H})\\ \bm{\chi}^{\ddagger} \odot \D {\bf b}}\right\|_{\infty}= \left\|\bmatrix{\bm{\Phi}_{\mathbb{S}} & \bf 0\\ \bf 0 & I_{\bm l}}  \bmatrix{\vec_{\mathbb{S}}(\Psi^{\ddagger} \odot\D \mathbf{H})\\ \bm{\chi}^{\ddagger} \odot \D {\bf b}}\right\|_{\infty}= \left\| \bmatrix{\vec_{\mathbb{S}}(\Psi^{\ddagger} \odot\D \mathbf{H})\\ \bm{\chi}^{\ddagger} \odot \D {\bf b}} \right\|_{\infty}.
\end{equation}
By considering $\tau=\gamma=\infty$ on \eqref{def:SCN} and using \eqref{eq:513}, \eqref{th1:eq2} and \eqref{th1:eq1}, we obtain
 \begin{align*}
         \mathfrak{K}^{\infty,\mathbb{S}}_{\vp}(\mathbf{H}, \b; \mathrm{L})&= \sup_{\substack{\left\|\vec\left(\Psi^{\ddagger} \odot \Delta \mathbf{H},\, \be^{\ddagger} \odot \D \b\right)\right\|_{\infty}\neq 0\\ \D A\in\, \mathbb{S}_1,\, \D D\in\, \mathbb{S}_2,\, \D E \in\, \mathbb{S}_3} } \frac{\left\|\Theta_{\xi^{\ddagger}_{\mathrm{L}}}\mathrm{L} 
 \mathfrak{B}^{-1} \bmatrix{ \mathcal{G} &  -I_{\bm l}} \bmatrix{\bm{\Phi}_{\mathbb{S}} & \bf 0\\ \bf 0 & I_{\bm l}} \bmatrix{\vec_{\mathbb{S}}( \D \mathbf{H})\\  \D {\bf b}}\right\|_{\infty}}{\left\|\vec\left(\Psi^{\ddagger} \odot \Delta \mathbf{H}, \be^{\ddagger} \odot \D \b\right)\right\|_{\infty}}\\
 &= \sup_{\footnotesize\substack{\left\|\bmatrix{\vec_{\mathbb{S}}(\Psi^{\ddagger} \odot \D \mathbf{H})\\ \bm{\chi}^{\ddagger} \odot \D {\bf b}}\right\|_{\infty}\neq 0\\ \D A\in\, \mathbb{S}_1,\, \D D\in\, \mathbb{S}_2,\, \D E \in\, \mathbb{S}_3} } \frac{\left\|\Theta_{\xi^{\ddagger}_{\mathrm{L}}}\mathrm{L} 
 \mathfrak{B}^{-1} \bmatrix{ \mathcal{G} &  -I_{\bm l}}  \bmatrix{\bm{\Phi}_{\mathbb{S}} & \bf 0\\ \bf 0 & I_{\bm l}} \bmatrix{\Theta_{\vec_{\mathbb{S}}(\Psi)}& \bf 0\\ \bf 0 & \Theta_{\be}}\bmatrix{\vec_{\mathbb{S}}(\Psi^{\ddagger} \odot\D \mathbf{H})\\ \bm{\chi}^{\ddagger} \odot \D {\bf b}}\right\|_{\infty}}{\left\|\bmatrix{\vec_{\mathbb{S}}(\Psi^{\ddagger} \odot \D \mathbf{H})\\ \bm{\chi}^{\ddagger} \odot \D {\bf b}}\right\|_{\infty}}\\
 &=\left\|\Theta_{\xi^{\ddagger}_{\mathrm{L}}}\mathrm{L} 
 \mathfrak{B}^{-1} \bmatrix{ \mathcal{G} &  -I_{\bm l}} \bmatrix{\bm{\Phi}_{\mathbb{S}} & \bf 0\\ \bf 0 & I_{\bm l}}\bmatrix{\Theta_{\vec_{\mathbb{S}}(\Psi)}& \bf 0\\ \bf 0 & \Theta_{\be}} \right\|_{\infty}\\
 &=\left\||\Theta_{\xi^{\ddagger}_{\mathrm{L}}}| \left|\mathrm{L}\mathfrak{B}^{-1} \bmatrix{ \mathcal{G} &  -I_{\bm l}} \bmatrix{\bm{\Phi}_{\mathbb{S}} & \bf 0\\ \bf 0& I_{\bm l}}\right| \bmatrix{\vec_{\mathbb{S}}(|\Psi|)\\  |\be|}\right\|_{\infty}.
    \end{align*}

    The rest of the proof follows by considering $\Psi=\mathbf{H},$ $\be ={\bf b},$ $\xi_{\mathrm{L}}=\|\mathrm{L}[\bm{x}^{\bm\top}, \bm{y}^{\bm\top}, \bm{z}^{\bm\top}]^{\bm\top}\|_{\infty}{\bf{1}}_{\bm l}$ (or $\xi_{\mathrm{L}}= \mathrm{L}[\bm{x}^{\bm\top}, \bm{y}^{\bm\top}, \bm{z}^{\bm\top}]^{\bm\top}$). 
$\blacksquare $

 
 \section{Deduction of Partial CNs for the EILS problem}\label{sec:EILS}
The EILS problem is an extension of the famous linear least squares problem, having linear constraints on unknown parameters. It  can be expressed as follows:
\begin{equation}\label{EILS}
    \min_{\y\in \R^m}(b-M\y)^{\bm\top}\mathbb{J}~(b-M\y)\quad \text{subjected to}\quad C\y=d,
\end{equation}
where $M\in \R^{n\times m} (n\geq m),$ $C\in \R^{p\times m},$ $b\in \R^n$ $d\in \R^p$ and the signature matrix $ \mathbb{J}$ given by
\begin{equation}
    \mathbb{J}=\bmatrix{I_{n_1} & \bf 0\\ \bf 0 & -I_{n_2}}, ~ n_1+n_2=n.
\end{equation}
When $\rank(C)=p$ and $\y^{\bm\top}(M^{\bm\top}\mathbb{J}M)\y>0$ for all nonzero $\y\in \text{null}(C),$ the EILS problem \eqref{EILS} has a unique solution. 
The solution of the EILS \eqref{EILS} problem also satisfies the following the augmented system \cite{ILSE2022}:
\begin{equation}\label{EILS2}
   \widehat{\mathfrak{B}}\bmatrix{\lambda \\ \x \\ \y}:= \bmatrix{\bf 0 & \bf 0 & C\\ \bf 0 & \mathbb{J} & M\\ C^{\bm\top} & M^{\bm\top}& \bf 0 }\bmatrix{\lambda \\ \x \\ \y}=\bmatrix{d\\b\\ \bf 0},
\end{equation}
where $\x=\mathbb{J}{\bm r},$ ${\bm r}=b-M\y$ and $\lambda=(CC^{\bm\top})^{-1}CM^{\bm\top}\mathbb{J}{\bm r}$ is the vector of Lagrange multipliers \cite{EILSAlgorithms}. 
Note that the system in \eqref{EILS2} can be equivalently transformed into 
\begin{align}\label{EILS3}
    \bmatrix{\mathbb{J} & M & \bf 0 \\ M^{\bm\top} & \bf 0 & C^{\bm\top}\\ \bf 0 & C & \bf 0}\bmatrix{\x \\ \y \\ \lambda}=\bmatrix{ b\\\bf 0\\ d }=:\widehat{\bf d}.
\end{align}
Observe that, the above system is in the form of DSPP \eqref{eq1:SPP} with $A=\mathbb{J},$ $B=M^{\bm\top},$ ${\bf b}=[b^{\bm\top}, {\bf 0}, d^{\bm\top}]^{\bm\top}$ and $\z=\lambda.$
Therefore, the task of assessing the conditioning of the EILS problem \eqref{EILS} can be achieved by determining the CNs for the solution $\y$ of DSPP \eqref{EILS3}.

Generally, the signature matrix $\mathbb{J}$ has no perturbation and as $D={\bf 0},$ $E={\bf 0}$ and ${\bf b}_2={\bf 0}$, we consider $\D A=\0$, $\D D=\0,$ $\D E=0$  and  $\D {\bf b}_2={\bf 0}$ in \eqref{eq:pert}. 
Then, the perturbation expression in \eqref{eq27} reduces to
\begin{eqnarray}\label{eq45}
  \bmatrix{\D\bm{x}\\ \D\bm{z}\\ \D\bm{y}}= -\mathfrak{B}^{-1} \bmatrix{\widehat{\mathcal{G}} &  -I_{n+p}}\bmatrix{\vec(\D B^{\T})\\ \vec(\D C)\\ \D {\bf b}_1 \\ \D {\bf b}_3} + \mathcal{O}(\bm{\epsilon}^2),
\end{eqnarray}
 where \begin{equation}
    \widehat{\mathcal{G}}= \bmatrix{ \bm{y}^{\bm\top}\otimes I_n& \bf 0\\  I_m\otimes \bm{x}^{\bm\top} & I_m\otimes \bm{z}^{\bm\top} \\ \bf 0& \bm{y}^{\bm\top}\otimes I_p }\in \R^{\bm{l}\times \hat{\bm s}},
 \end{equation}
 $\hat{\bm s}=m(n+p).$

Let $\widehat{\mathbf{H}}=(B^{\T},C),$ $\D\widehat{\mathbf{H}}=(\D B^{\T}, \D C),$ $\widehat{\bf b}=[{\bf b}_1^{\bm\top}, {\bf b}_3^{\bm\top} ]^{\bm\top}$ and $\D \widehat{\bf b}=[\D {\bf b}_1^{\bm\top}, \D{\bf b}_3^{\bm\top} ]^{\bm\top}$.
We define the following mapping:
\begin{align}\label{eq:map2}
\nonumber \widetilde{\vp}:\,&   \R^{m\times n}\times  \R^{p\times m}  \times \R^{n+p  } \rightarrow \R^k\\
&~~\widetilde{\vp}(\widehat{\mathbf{H}}, \widehat{\mathbf{b}})=\mathrm{L}\bmatrix{\x\\ \y \\ \z}=\mathrm{L}\mathfrak{B}^{-1} \bmatrix{{\bf b}_1\\ \0\\  {\bf b}_3},
\end{align}
where $\mathrm{L}\in \R^{k\times \bm l}.$ Using a similar method to the Theorem \ref{th1}, we have the following result.
 \begin{theorem}\label{Th41}
 Assume that  $[\bm{x}^{\bm\top}, \bm{y}^{\bm\top}, \bm{z}^{\bm\top}]^{\bm\top}$ is the unique solution of the DSPP \eqref{eq1:SPP} with $D=\0,$ $E=\0$,  ${\bf b}_2={\bf 0}$ and $\mathrm{L}\in \R^{k\times \bm{l}}.$ Then the partial unified CN of $[\bm{x}^{\bm\top}, \bm{y}^{\bm\top}, \bm{z}^{\bm\top}]^{\bm\top}$ w.r.t. $\mathrm{L}$ is given by
    \begin{align*}
         \mathfrak{K}_{\widetilde{\vp} }(\widehat{\mathbf{H}}, \b; \mathrm{L})&:=\lim_{\bm{\epsilon}\rightarrow 0} \sup_{0<\left\|\vec\left(\widehat{\Psi}^{\ddagger} \odot \Delta \widehat{\mathbf{H}},\, \widehat{\be}^{\ddagger} \odot \D \b\right)\right\|_{\tau}\leq \bm{\epsilon}} \frac{\left\|\xi^{\ddagger}_{\mathrm{L}}\odot \left( \widetilde{\bm{\varphi}}(\widehat{\mathbf{H}}+\D \widehat{\mathbf{H}}, \widehat{\mathbf{b}}+ \D \widehat{\mathbf{b}})-\widetilde{\vp}(\widehat{\mathbf{H}}, \widehat{\b})\right) \right\|_\gamma}{\left\|\vec\left(\widehat{\Psi}^{\ddagger} \odot \Delta \widehat{\mathbf{H}}, \widehat{\be}^{\ddagger} \odot \D \widehat{\b}\right)\right\|_{\tau}},\\
         &\, = \left\|\Theta_{\xi^{\ddagger}_{\mathrm{L}}}\mathrm{L}\mathfrak{B}^{-1}\bmatrix{ \widehat{\mathcal{G}} &  -I_{n+p}} \bmatrix{\Theta_{\vec(\widehat{\Psi})}& \bf 0\\ \bf 0 & \Theta_{\widehat{\be}}}\right\|_{\tau,\gamma},
             \end{align*}
    where $\widehat{\Psi}= (\Psi_{B^{\T}},\Psi_C),$ $\Psi_{B^{\T}}\in \R^{n\times m}, \Psi_C\in \R^{p\times m}$ and $\widehat{\be}\in \R^{n+p}.$ 
 \end{theorem}

\begin{remark}
Taking $\mathrm{L}=[{\0}_{k\times n} ~~ \mathrm{L}_1~ ~{\0}_{k\times p}],$ $\mathrm{L}_1\in \R^{k\times m},$ $A=\mathbb{J}$  in Theorem \ref{Th41}, and since $\mathfrak{B} = \Sigma \widehat{\mathfrak{B}} \Sigma^{-1},$ where
\begin{equation}
    \Sigma=\bmatrix{\0&I_n &\0  \\ \0 & \0 &I_m\\ I_p &\0& \0},
\end{equation}
using the formula for the inverse of $\widehat{\mathfrak{B}}$ given in \cite{EILS2010}, we obtain
\begin{align}
       \mathfrak{K}_{\widetilde{\vp} }(\widehat{\mathbf{H}}, \b; \mathrm{L})=\left\|\Theta_{\xi^{\ddagger}_{\mathrm{L}}} \mathrm{L}_1 \bmatrix{\Xi& \Lambda&  -(\mathcal{Q}\mathcal{M}\mathcal{Q})^{\dagger}MJ & \mathcal{B}_M } \bmatrix{\Theta_{\vec(\widehat{\Psi})}& \bf 0\\ \bf 0 & \Theta_{\widehat{\be}}} \right\|_{\tau,\gamma},
\end{align}
where $\Xi=\y^{\bm\top} \otimes (\mathcal{Q}\mathcal{M}\mathcal{Q})^{\dagger}M\mathbb{J}- (\mathcal{Q}\mathcal{M}\mathcal{Q})^{\dagger}\otimes \x^{\bm\top},$ $\Lambda=\y^{\bm\top}\otimes \mathcal{B}_{M}-(\mathcal{Q}\mathcal{M}\mathcal{Q})^{\dagger}\otimes \z^{\bm\top},$ $\mathcal{M}=M\mathbb{J}M^{\bm\top},$ $\mathcal{Q}=I_m-C^{\dagger}C$ and $\mathcal{B}_M=(I_m-\mathcal{Q}\mathcal{M}\mathcal{Q})^{\dagger}.$ Note that the partial CN expression is the same as derived in \cite{ILSE2022}.
    
\end{remark}

\section{Numerical experiments}\label{sec:Numerical}
In this part, we present some numerical examples to verify the reliability of the derived partial NCN, MCN, and CCN and their upper bounds for the DSPP. Additionally, we demonstrate their effectiveness in providing tight upper bounds for the relative forward error in solving the DSPP.

We construct the entrywise perturbation as follows:
\begin{align*}
   & \D A=10^{-s}\cdot \mathtt{randn}(n,n)\odot A,~  \D B=10^{-s}\cdot \mathtt{randn}(m,n)\odot B,~ \D C=10^{-s}\cdot \mathtt{randn}(p,m)\odot C,\\
    & \D D=10^{-s}\cdot \mathtt{randn}(m, m)\odot D,~  \D E=10^{-s}\cdot \mathtt{randn}(p,p)\odot E,~ \D{\bf b}=10^{-s}\cdot \mathtt{randn}(\bm{l}, 1)\odot{\bf b},
\end{align*}
where $\mathtt{randn}(m,n)$ denotes $m\times n$ the random matrices generated by MATLAB command $\mathtt{randn}$. Let $\bm{w}= [\x^{\T},\y^{\T},\z^{\bm\top}]^{\T} $ and $\widetilde{\bm{w}}=[\widetilde{\x}^{\T}, \widetilde{\y}^{\T}, \widetilde{\z}^{\T}]^{\T}$ be the unique solutions of the original DSPP and the perturbed DSPP, respectively. To estimate an upper bound for the forward error in the solution, the normwise, mixed and componentwise relative forward errors  are defined as follows:
\begin{align*}
		\bm{r}_k=\frac{\|\mathrm{L}\widetilde{\bm{w}}-\mathrm{L}\bm{w}\|_2}{\|\mathrm{L}\bm{w}\|_2}, \quad  \bm{r}_m=\frac{\|\mathrm{L}\widetilde{\bm{w}}-\mathrm{L}\bm{w}\|_{\infty}}{\|\mathrm{L}\bm{w}\|_{\infty}}  \quad \text{and} ~~\bm{r}_c=\left\|\frac{\mathrm{L}\widetilde{\bm{w}}-\mathrm{L}\bm{w}}{\mathrm{L}\bm{w}}\right\|_{\infty},
	\end{align*}
 respectively.

{Similar to the relation in \eqref{eq:forward} and from Definition \ref{def:PartialCN}, the forward error can be easily bounded as follows:
\begin{eqnarray}\label{forward2}
       \left\|\xi^{\ddagger}_{\mathrm{L}}\odot \left( \bm{\varphi}(\mathbf{H}+\D \mathbf{H}, \mathbf{b}+ \D \mathbf{b})-\vp(\mathbf{H}, \b)\right) \right\|_\gamma\leq \bm{\epsilon} \,  \mathfrak{K}_{\vp}(\mathbf{H}, \b; \mathrm{L})
       ~+~\mathcal{O}(\bm{\epsilon}^2).
    \end{eqnarray}}
     Define the following quantities:
\begin{equation}
    \bm{\epsilon}_1=\frac{\left\|\bmatrix{\D \mathfrak{B}& \D \mathbf{b}}\right\|_F}{\left\|\bmatrix{\mathfrak{B} & \mathbf{b}}\right\|_F}~~ \text{and}~~\bm{\epsilon}_2=\min\left\{\bm{\epsilon}: \left|\D \mathfrak{B}\right|\leq \bm{\epsilon}\left|\mathfrak{B} \right|,~ \left|\D \mathbf{b}\right|\leq \bm{\epsilon}\left|\mathbf{b}\right|\right\}.
\end{equation}
{Then using the bound in \eqref{forward2} and using the setup of NCN, MCN and CCN}, $\bm{\epsilon}_1\mathfrak{K}^{(2)}_{\vp}(\mathbf{H}, \b; \mathrm{L}),$ $\bm{\epsilon}_2 \mathfrak{K}^{\infty}_{mix, \vp}(\mathbf{H}, \b; \mathrm{L})$ and $\bm{\epsilon}_2 \mathfrak{K}^{\infty}_{com, \vp}(\mathbf{H}, \b; \mathrm{L})$ can be employed to estimate the upper bound for relative forward errors $\bm{r}_k,$ $\bm{r}_m$ and $\bm{r}_c,$ respectively. We select the matrix \(\mathrm{L}\) as $$\mathrm{L}_0=I_{\bm l}, ~ \mathrm{L}_n=\bmatrix{ I_n & \mathbf{0}_{n\times (m+p)}},~ \mathrm{L}_m=\bmatrix{{\0}_{m\times n}& I_m & \mathbf{0}_{m\times p}} ~\text{and} ~ \mathrm{L}_p=\bmatrix{{\0}_{p\times (n+m)}& I_p }$$ to obtain the CNs for $\bm{w}$, \(\x\), $\y$ and $\z$, respectively.

\vspace{2mm}
\begin{exam}\label{exam1}
     We consider the DSPP \eqref{eq1:SPP} taken from \cite{HuangNA} with
		\begin{eqnarray*}
			& A= \bmatrix{I_q\otimes J+J\otimes I_q &\bf 0\\ \bf 0&I_q\otimes J+ J\otimes I_q}\in \R^{2q^2\times 2 q^2}, ~~  \\
	  &B=\bmatrix{I_q\otimes Z& Z\otimes I_q}\in \R^{q^2\times 2q^2},~ C= Y\otimes Z\in \R^{q^2\times q^2},
		\end{eqnarray*}
		where $J=\frac{1}{(q+1)^2}\, \mathtt{tridiag}(-1,2,-1)\in \R^{q\times q},$ $ Z=\frac{1}{q+1}\,  \mathtt{tridiag}(0,1,-1)\in \R^{q\times q}$  and  $Y=\diag(1, q+1, \ldots, q^2-q+ 1)\in \R^{q\times q}.$ Further, we take $D=I_m$ and $E=I_p.$ The notation $ \mathtt{tridiag}(d_1,d_2,d_3)\in \R^{q\times q}$ denotes the 
 tridiagonal matrix with subdiagonal entries $d_1$, diagonal entries $d_2$ and superdiagonal entries $d_3$. Here, the dimension of the coefficient matrix $\mathfrak{B}$ of the DSPP is ${\bm l}=4q^2.$   We take ${\bf b}=\tt{randn}(n,1)\in \R^{\bm l}$. Further, we consider $\Psi=\|\mathfrak{B}\|_F$ and ${\bm \chi}=\|{\bf b}\|_2.$ We use Theorems \ref{th32} and \ref{th33} to compute partial  CNs and Corollaries \ref{cor31} and \ref{coro32} for their upper bounds. {In Table \ref{tab:new}, we listed the computed values of $\bm{\epsilon}_1$ and $\bm{\epsilon}_2$. The numerical results of the different relative forward errors and their upper bounds for different choices for $\mathrm{L}$ and $q=4:2:16$ are presented in Tables \ref{tab1}-\ref{tab4}.} 
  \begin{table}[h]
\caption{Values of $\bm{\epsilon}_1$ and $\bm{\epsilon}_2$ for different $q.$}\label{tab:new}%
\begin{tabular}{@{}ccc@{}}
\toprule
 $q$   & $\bm{\epsilon}_1$ & $\bm{\epsilon}_2$\\
\midrule
  $4$ & $1.1234E-08$ & $3.2043E-08$       \\
  \midrule
  $6$ &$1.0007E-08$&	$3.4305E-08$ \\
  \midrule
  $8$ & $9.2718E-09$	&$2.9944E-08$\\
  \midrule
  $10$ & $1.0404E-08$	& $4.0272E-08$\\
  \midrule
  $12$ &$1.0384E-08$&	$3.8205E-08$ \\
\botrule
\end{tabular}
\end{table}
 \begin{sidewaystable}
	\centering
		\caption{Comparison of relative forward errors and their upper bounds obtained using the NCN, MCN, and CCN, and their upper bounds for   \(\mathrm{L} = \mathrm{L}_0\) for Example \ref{exam1}.}
			\label{tab1}
				\begin{tabular}{@{}cccccccccc@{}}
					\toprule
			$q$	& $\bm{r}_k$ &$ {U_{\bm{r}_k}}$	& $ {\widehat{U}_{\bm{r}_k}}$ & $\bm{r}_m$ &$ {U_{\bm{r}_m}}$& $ {\widehat{U}_{\bm{r}_m}}$ & $\bm{r}_c$ &  $ {U_{\bm{r}_c}}$& $ {\widehat{U}_{\bm{r}_c}}$ \\
  \midrule 
   4 &$1.1882E-08$&$5.1234E-06$&	$6.3888E-06$&	$1.9907E-08$&	$3.5552E-07$&	$3.8637E-07$&	$2.1532E-07$&	$6.9184E-06$&	$8.6300E-06$\\[1ex]
   $6$ &$2.0426E-08$&	$8.2886E-06$&	$1.6984E-05$&	$2.3603E-08$&	$3.9780E-07$&	$4.0885E-07$&	$3.7248E-07$&	$1.2135E-05$&	$1.2490E-05$\\[1ex]
     $8$ &$3.1951E-08$&	$1.9154E-05$&	$3.6586E-05$&	$4.1045E-08$&	$4.9025E-07$&	$4.9808E-07$&	$8.3715E-06$&	$1.3746E-04$&	$1.3964E-04$\\	[1ex]
    $10$ &$2.9187E-08$&	$5.5170E-05$&	$1.1519E-04$&	$4.2771E-08$&	$1.2319E-06$&	$1.2706E-06$&	$1.8721E-07$&	$1.9213E-05$&	$2.1515E-05$\\[1ex]
     $12$ &$2.3289E-08$&	$9.4616E-05$&	$1.7034E-04$&	$2.9289E-08$&	$1.0722E-06$&	$1.0798E-06$&	$4.5669E-07$&	$3.5850E-05$&	$3.8136E-05$\\[1ex]
   $14$ & $3.9431E-08$&	$1.6571E-04$&	$3.3889E-04$&	$3.9739E-08$&	$1.2469E-06$&	$1.2747E-06$&	$6.0883E-07$&	$3.1225E-05$&$	3.5558E-05$\\[1ex]
  $16$ &$3.9038E-08$&	$2.7496E-04$&	$5.6909E-04$&	$3.9932E-08$&	$1.8001E-06$&	$1.8356E-06$&	$4.2789E-06$&	$1.7386E-04$&$	1.8883E-04$\\[1ex]
   \bottomrule   
    \end{tabular} 
    \end{sidewaystable}
     \begin{sidewaystable}
			\centering
		\caption{Comparison of relative forward errors and their upper bounds obtained using the NCN, MCN, and CCN, and their upper bounds    \(\mathrm{L} = \mathrm{L}_n\) for Example \ref{exam1}.}
			\label{tab2}
				\begin{tabular}{cccccccccc}
					\toprule
			$q$	& $\bm{r}_k$ &$ {U_{\bm{r}_k}}$	& $ {\widehat{U}_{\bm{r}_k}}$ & $\bm{r}_m$ &$ {U_{\bm{r}_m}}$& $ {\widehat{U}_{\bm{r}_m}}$ & $\bm{r}_c$ &  $ {U_{\bm{r}_c}}$& $ {\widehat{U}_{\bm{r}_c}}$ \\
  \midrule 
   $4$ &$2.3210E-08$&	$1.7899E-06$&	$2.2325E-06$&	$4.4418E-08$	&$3.9328E-07$	&$3.9952E-07$&	$6.1112E-08$&	$7.7010E-07$	&$7.9196E-07$\\[1ex]
   $6$ &$2.2133E-08$&	$2.6980E-06$&	$7.6732E-06$&	$4.2425E-08$&	$8.8875E-07$&	$9.5134E-07$&	$5.3918E-08$&	$5.8823E-06$&	$7.3855E-06$\\[1ex]
     $8$ &$3.8470E-08$&	$4.8135E-06$&	$1.5701E-05$&	$8.7676E-08$&	$1.4710E-06$	&$1.5322E-06$	&$5.7409E-07$	&$5.8895E-05$&	$7.0437E-05$\\	[1ex]
    $10$ &$2.9894E-08$&	$9.6923E-06$&	$3.1204E-05$&	$5.9374E-08$	&$2.0857E-06$	&$2.1065E-06$&	$2.7350E-07$&	$9.5571E-05$	&$1.0598E-04$\\[1ex]
     $12$& $3.3113E-08$&$3.6464E-05$&	$4.9737E-05$&	$5.6822E-08$&	$3.0914E-06$&	$3.0928E-06$&	$4.5080E-07$&	$9.3747E-06$&	$9.4743E-06$\\[1ex]
   $14$ & $8.0906E-08$&	$2.4824E-05$	&$7.6132E-05$	&$1.5910E-07$&	$3.3988E-06$&	$3.4220E-06$&	$8.4585E-07$&	$2.5939E-05$&	$2.9725E-05$\\[1ex]
  $16$ &$4.4663E-08$&	$3.5262E-05$&	$8.8133E-05$	&$8.0314E-08$&	$4.7189E-06$&	$4.7301E-06$	&$4.5566E-07$	&$3.7542E-05$	&$3.9779E-05$\\[1ex]
   \bottomrule
    \end{tabular}
    %
    \end{sidewaystable}
    
    \begin{sidewaystable}
	\footnotesize		\centering
		\caption{Comparison of relative forward errors and their upper bounds obtained using the NCN, MCN, and CCN, and their upper bounds for   \(\mathrm{L} = \mathrm{L}_m\) for Example \ref{exam1}.}
			\label{tab3}
				\begin{tabular}{cccccccccc}
					\toprule
			$q$	& $\bm{r}_k$ &$ {U_{\bm{r}_k}}$	& $ {\widehat{U}_{\bm{r}_k}}$ & $\bm{r}_m$ &$ {U_{\bm{r}_m}}$& $ {\widehat{U}_{\bm{r}_m}}$ & $\bm{r}_c$ &  $ {U_{\bm{r}_c}}$& $ {\widehat{U}_{\bm{r}_c}}$  \\
  \midrule 
   $4$ &$1.4640E-08$&	$4.1828E-06$&	$8.1301E-06$	&$1.7030E-08$&	$3.2527E-07$&	$3.6699E-07$&	$6.3575E-08$	&$1.6374E-06$&	$1.7734E-06$\\[1ex]
   $6$& $3.8057E-08$&	$3.4032E-05$&	$5.4748E-05$&	$3.9292E-08$	&$8.5355E-07$	&$1.0120E-06$	&$6.8320E-08	$&$2.2698E-06$	&$2.3143E-06$\\[1ex]
     $8$ &$6.2425E-08$	&$1.1653E-04$&	$1.7568E-04$&	$5.3397E-08$&	$1.9931E-06$&	$2.1774E-06$&	$3.0610E-07$	&$1.3012E-05$	&$1.3193E-05$\\	[1ex]
    $10$ &$4.7836E-08$&	$1.5526E-04$&	$3.3655E-04$&	$5.1330E-08$&	$1.8322E-06$&	$2.2332E-06$&	$5.0917E-07$&	$2.8394E-05$&	$2.8957E-05$\\[1ex]
     $12$&$ 3.7645E-08$ &	$9.3405E-05$&	$4.7395E-04$	&$3.6702E-08$	&$1.1129E-06$&	$1.4824E-06$	&$7.7646E-08$	&$7.9413E-06$&	$8.8422E-06$\\[1ex]
   $14$ & $7.7970E-08$&	$1.6058E-04$&	$9.2172E-04$&	$5.7666E-08$&	$1.7296E-06$&	$2.1704E-06$&	$2.9851E-06$&	$7.6483E-05$	&$7.9562E-05$\\[1ex]
  $16$ &$4.7259E-08$&	$4.3435E-04$&	$1.5265E-03$&	$6.6721E-08$&	$2.2484E-06$	& $2.8853E-06$&	$2.9895E-07$&	$5.8651E-05$&	$6.0528E-05$\\[1ex]
   \bottomrule

    \end{tabular}
    \end{sidewaystable}

     \begin{sidewaystable}
		\centering
		\caption{Comparison of relative forward errors and their upper bounds obtained using the NCN, MCN, and CCN, and their upper bounds for   \(\mathrm{L} = \mathrm{L}_p\) for Example \ref{exam1}.}
			\label{tab4}
				\begin{tabular}{@{}cccccccccc@{}}
					\toprule
			$q$	& $\bm{r}_k$ &$ {U_{\bm{r}_k}}$	& $ {\widehat{U}_{\bm{r}_k}}$ & $\bm{r}_m$ &$ {U_{\bm{r}_m}}$& $ {\widehat{U}_{\bm{r}_m}}$ & $\bm{r}_c$ &  $ {U_{\bm{r}_c}}$& $ {\widehat{U}_{\bm{r}_c}}$ \\
  \midrule 
   $4$ &$7.9587E-09$	& $3.2626E-06$ &	$5.0845E-06$	& $9.2868E-09$	&$2.7450E-07$&	$2.8678E-07$&	$6.4428E-08$&	$5.9782E-07$	&$6.2833E-07$\\[1ex]
   $6$&$1.3913E-08$&$	9.6919E-06$&	$1.1609E-05$&	$1.8364E-08$&	$4.4130E-07$	&$4.4401E-07$&	$4.1485E-08$&	$1.1621E-06$	&$1.1751E-06$\\[1ex]
     $8$ &$2.1543E-08$&	$2.5173E-05$&	$4.6780E-05$	&$2.7719E-08$&	$7.9552E-07$&	$8.4231E-07$&	$4.9306E-08$&	$4.2242E-06$&	$4.4848E-06$\\	[1ex]
    $10$ &$3.2132E-08$&	$5.7696E-05$&	$1.1250E-04$&	$4.3437E-08$&	$1.0618E-06$	&$1.1161E-06$&	$8.3724E-08$&	$4.1932E-06$&	$4.3258E-06$\\[1ex]
     $12$&$1.6735E-08$&	$8.7652E-05$&	$1.4248E-04$&	$3.1936E-08$&	$1.0011E-06$&	$1.0138E-06$&	$6.2623E-08$&	$2.5824E-06$	&$2.5860E-06$\\[1ex]
   $14$ & $2.0065E-08$&	$1.5497E-04$&	$1.9181E-04$&	$2.2200E-08$&	$1.0421E-06$&	$1.0457E-06$&	$7.7378E-08$&	$4.1703E-06$&	$4.1741E-06$\\[1ex]
  $16$ &$1.0902E-08$&	$2.4487E-04$	&$3.5342E-04$	&$1.5817E-08$	&$1.6369E-06$&$	1.6456E-06$&	$1.0229E-07$&	$3.4527E-06$&	$3.4614E-06$\\[1ex]
   \bottomrule

    \end{tabular}
    \end{sidewaystable}
  For simplicity,  {in Tables \ref{tab1}-\ref{tab4},} we denote:
 { \begin{align*}
     & {U_{\bm{r}_k}}=\bm{\epsilon}_1\mathfrak{K}^{(2)}_{\vp}(\mathbf{H}, \b; \mathrm{L}), ~ {\widehat{U}_{\bm{r}_k}}=\bm{\epsilon}_1\mathfrak{K}^{(2),u}_{\vp}(\mathbf{H}, \b; \mathrm{L}),~
         {U_{\bm{r}_m}}=\bm{\epsilon}_2\mathfrak{K}^{\infty}_{mix,  \vp}(\mathbf{H}, \b; \mathrm{L}),\\
       &  {\widehat{U}_{\bm{r}_m}}=\bm{\epsilon}_2\mathfrak{K}^{\infty,u}_{mix,  \vp}(\mathbf{H}, \b; \mathrm{L}),~
         {U_{\bm{r}_c}}=\bm{\epsilon}_2\mathfrak{K}^{\infty}_{com,  \vp}(\mathbf{H}, \b; \mathrm{L}),~
         {\widehat{U}_{\bm{r}_c}}=\bm{\epsilon}_2\mathfrak{K}^{\infty,u}_{com,  \vp}(\mathbf{H}, \b; \mathrm{L}).
    \end{align*}}
     {By \eqref{eq:forward}, the quantities $U_{\bm{r}_k}$ and $\widehat{U}_{\bm{r}_k}$ are the upper bounds of the relative forward error $\bm{r}_k$, $U_{\bm{r}_m}$ and $\widehat{U}_{\bm{r}_m}$ are the upper bounds of the relative forward error $\bm{r}_m$, and $U_{\bm{r}_c}$ and $\widehat{U}_{\bm{r}_c}$ are the upper bounds of the relative forward error $\bm{r}_c.$ These bounds are obtained using the NCN, MCN, CCN, and their upper bounds, respectively.}

     {As by multiplying the same perturbation magnitude \( \bm{\epsilon}_1 \) with the NCN and its upper bound we obtain \( U_{\bm{r}_k} \) and \( \widehat{U}_{\bm{r}_k} \), respectively, and as shown in Tables~2--5, they consistently remain of the same order across all values of \( q \). Moreover, \( U_{\bm{r}_k} \) being slightly smaller than \( \widehat{U}_{\bm{r}_k} \). This demonstrates that the upper bound of the NCN provides a sharp estimate of the NCN. By a similar argument, we found that the upper bounds of the MCN and CCN are also sharp estimates of the MCN and CCN.
}

   Additionally, it is observed that both the MCN and CCN, along with their upper bounds, are at most two orders of magnitude larger than the actual relative forward errors, offering more accurate estimates compared to the NCN and its upper bounds. These numerical results highlight the effectiveness of the proposed CNs and their corresponding upper bounds.

   {To evaluate the robustness of the proposed framework for computing the partial NCN, MCN and CCN of the DSPP with respect to small perturbations in the input data, we consider a perturbed coefficient matrix defined as $\widehat{\mathfrak{B}} = \mathfrak{B} + 10^{-4} \Delta \mathfrak{B}$, where the block matrices of $\Delta \mathfrak{B}$ are specified as follows:
   \begin{align*}
      & \D A=e_p \,{\tt std}(A) \,{\tt{randn}}(n,n),~\D B=e_p \,{\tt std}(B) \,{\tt{randn}}(m,n),~\D C=e_p \,{\tt std}(C) \,{\tt{randn}}(p,m),\\
       &\D D=e_p \,{\tt std}(D) \,{\tt{randn}}(m,m),~\text{and}~~\D E=e_p \,{\tt std}(E) \,{\tt{randn}}(p,p),
   \end{align*}
where   ${\tt std}(A)$ denotes the standard deviation of $A.$ For the perturbed coefficient matrix, we consider the perturbed DSPP $\widehat{\mathfrak{B}}\widehat{\bm w}=\b.$
 \begin{table}[h]
\caption{Partial NCN, MCN and CCN for the original and perturbed DSPP by varying $e_p$ for $q=4.$ }\label{tab:new2}%
 \begin{tabular}{@{}cccc|ccc@{}}
\toprule
&\multicolumn{3}{c}{Original DSPP}&\multicolumn{3}{|c}{Perturbed DSPP}\\
\midrule
 $e_p$   & $\mathfrak{K}^{(2)}_{\vp}(\mathbf{H}, \b; \mathrm{L})$ & $\mathfrak{K}^{\infty}_{mix,\vp}(\mathbf{H}, \b; \mathrm{L})$& $\mathfrak{K}^{\infty}_{com,\vp}(\mathbf{H}, \b; \mathrm{L})$&$\mathfrak{K}^{(2)}_{\vp}(\mathbf{H}, \b; \mathrm{L})$&$\mathfrak{K}^{\infty}_{mix,\vp}(\mathbf{H}, \b; \mathrm{L})$& $\mathfrak{K}^{\infty}_{com,\vp}(\mathbf{H}, \b; \mathrm{L})$\\
\midrule
  $0.10$ & $2.0808E+02$ & $8.3003E+00$& $1.2487E+01$&$2.0809E+02$&$8.3003E+00$& $1.2487E+01$       \\
  \midrule
  $0.15$ &$6.2669E+02$&	$1.4536E+01$ & $1.6732E+01$&$6.2669E+02$&	$1.4536E+01$ & $1.6732E+01$\\
  \midrule
  $0.20$ & $1.0318E+03$	&$2.0297E+01$ & $2.8747E+02$& $1.0319E+03$& $2.0297E+01$& $2.8742E+02$\\
\botrule
\end{tabular}
\end{table}
\begin{table}[h]
\caption{Partial NCN, MCN and CCN for the original and perturbed DSPP by varying $e_p$ for $q=6.$ }\label{tab:new3}%
 \begin{tabular}{@{}cccc|ccc@{}}
\toprule
&\multicolumn{3}{c}{Original DSPP}&\multicolumn{3}{|c}{Perturbed DSPP}\\
\midrule
 $e_p$   & $\mathfrak{K}^{(2)}_{\vp}(\mathbf{H}, \b; \mathrm{L})$ & $\mathfrak{K}^{\infty}_{mix,\vp}(\mathbf{H}, \b; \mathrm{L})$& $\mathfrak{K}^{\infty}_{com,\vp}(\mathbf{H}, \b; \mathrm{L})$&$\mathfrak{K}^{(2)}_{\vp}(\mathbf{H}, \b; \mathrm{L})$&$\mathfrak{K}^{\infty}_{mix,\vp}(\mathbf{H}, \b; \mathrm{L})$& $\mathfrak{K}^{\infty}_{com,\vp}(\mathbf{H}, \b; \mathrm{L})$\\
\midrule
  $0.10$ & $4.2282E+03$ & $3.3861E+01$& $1.8403E+03$&$4.2283E+03$&$3.3861E+01$& $1.8396E+03$       \\
  \midrule
  $0.15$ &$5.9136E+03$&	$4.8040E+01$ & $1.5365E+02$&$5.9136E+03$&	$4.8041E+01$ & $1.5366E+02$\\
  \midrule
  $0.20$ & $1.1828E+03$	&$1.1395E+01$ & $2.3067E+01$& $1.1829E+03$& $1.1395E+01$& $2.3067E+01$\\
\botrule
\end{tabular}
\end{table}

The computed partial NCN, MCN and CCN for both the original and perturbed DSPP are presented in Tables \ref{tab:new2} and \ref{tab:new3} for $e_p=0.10, 0.15, 0.20$. The perturbation introduced in the coefficient matrices is of the order $\mathcal{O}(10^{-3})$. From these tables, we observe that the partial NCN, MCN, and CCN remain nearly identical across all cases. This demonstrates that the proposed expressions for the partial CNs are robust and exhibit low sensitivity to small perturbations in the input data.
} 
\end{exam}

\begin{exam}\label{example2}
We consider the DSPP \eqref{eq1:SPP} with the block matrices given by
    \begin{align*}
 &A=[(2ZZ^{\T}+\Sigma_1) \oplus \Sigma_2\oplus \Sigma_3]\in \R^{n\times n}, ~ B=\bmatrix{N&-I_{2\tilde{q}}& I_{2\tilde{q}}}\in \R^{m\times n},~ \\
     &D={\tt{toeplitz}}(\bm d)\in \R^{m\times m}, ~C=M, ~E={\tt{toeplitz}}(\bm e)\in \R^{p\times p}, \\
     &~\text{where}~  Z=[z_{ij}]\in \R^{\hat{q}\times \hat{q}}~ \text{with}~ z_{ij}=e^{-2((i/3)^2+(j/3)^2)}, ~\Sigma_1=I_{\hat{q}}, 
 \end{align*}
 \vspace{-3mm}
and $\Sigma_k=\diag(d_j^{(k)})\in \R^{2\tilde{q}\times 2\tilde{q}}, ~k=2,3,$ are diagonal matrices with  
$$ d_{j}^{(2)}=\left\{\begin{array}{cc}
    1,& \text{for} ~1\leq j\leq \tilde{q},  \\
    10^{-5}(j-\tilde{q})^2, & ~\text{for}~ \tilde{q}+1\leq j\leq 2\tilde{q},
\end{array}\right.$$
$d^{(3)}_j=10^{-5}(j+\tilde{q})^2$ for $1\leq j\leq 2\tilde{q},$ where $\tilde{q}=q^2$ and $\hat{q}=q(q+1).$ Further, $N=\bmatrix{\widehat{N}\otimes I_q\\ I_q\otimes \widehat{N}}\in \R^{2\tilde{q}\times \hat{q}},$ $\widehat{N}={\tt tridiag}(0,2,-1)\in \R^{q\times(q+1)},$ $M=\bmatrix{\widehat{M}\otimes I_q & I_q\otimes \widehat{M}}\in \R^{\hat{q}\times2\tilde{q}},$ 
$$\widehat{M}=\bmatrix{q+1& -\frac{q-1}{q}& \frac{q-2}{q}&\ldots &\frac{(-1)^{q-1}}{q}\\
-\frac{q-1}{q}& 2q+1& \ddots&\ddots&\vdots\\
\frac{q-2}{q}&\ddots &\ddots&\ddots&\frac{q-2}{q}\\
\vdots&\ddots&\ddots&\ddots&-\frac{q-1}{q}\\\frac{(-1)^{q-1}}{q}&\vdots&\frac{q-2}{q}& -\frac{q-1}{q}& q^2+1\\0 &0&\ldots& 0&1
}\in \R^{(q+1)\times q},$$
$\bm{d}=\mathtt{randn}(m,1)\in \R^m, \bm{e}=\mathtt{randn}(p,1)\in \R^p.$
Here, the notation $[A\oplus B]$ means the block diagonal matrix with diagonal blocks $A$ and $B,$ and $\mathtt{toeplitz}(\bm{d})$ denotes the Toeplitz matrix generated by the vector $\bm{d}.$ Moreover, we get ${\bm{l}}=8q^2+2q.$ The vector ${\bf d}\in\R^{\bm l}$ is chosen as in Example \ref{exam1}.

\begin{sidewaystable}
		\centering
		\caption{Comparison of the partial NCN, MCN, and CCN with the corresponding structured partial NCN, MCN and CCN for Example \ref{example2}.}
			\label{tab5}
				\begin{tabular}{cccccccc}
					\toprule
			$\mathrm{L}$&$q$	 &$\mathfrak{K}^{(2)}_{\vp}(\mathbf{H}, \b; \mathrm{L})$	& $\mathfrak{K}^{(2),\mathbb{S}}_{\vp}(\mathbf{H}, \b; \mathrm{L})$ & $\mathfrak{K}^{\infty}_{mix,  \vp}(\mathbf{H}, \b; \mathrm{L})$& $\mathfrak{K}^{\infty,\mathbb{S}}_{mix,  \vp}(\mathbf{H}, \b; \mathrm{L})$ &  $\mathfrak{K}^{\infty}_{com,  \vp}(\mathbf{H}, \b; \mathrm{L})$ &  $\mathfrak{K}^{\infty,\mathbb{S}}_{com,  \vp}(\mathbf{H}, \b; \mathrm{L})$ \\
  \midrule 
\multirow{3}{*}{$\mathrm{L}_0$}  & $2$ &$9.0961E+03$&	$8.9842E+03$&	$1.1174E+02$&	$9.2003E+01$&	$1.0211E+03$	&$7.8013E+02$\\[1ex]
   &$3$ &$1.0996E+04$&	$8.6963E+03$&	$1.1735E+02$&	$9.1218E+01$&	$1.2372E+04$&	$9.5905E+03$	\\[1ex]
   &  $4$ &$2.6066E+04$&	$1.9983E+04$	&$3.5456E+02$&	$2.6052E+02$&	$2.5284E+05$	&$6.7124E+04$	\\	[1ex]
   \midrule
   \multirow{3}{*}{$\mathrm{L}_n$}  & $2$ &$9.1011E+03$&	$8.9891E+03$	&$1.1174E+02$	&$9.2003E+01$	&$2.9170E+02$&	$2.1175E+02$\\[1ex]
   &$3$ &	$1.0972E+04$	&$8.6781E+03$&	$1.1735E+02$&	$9.1218E+01$&	$1.2372E+04$&	$9.5905E+03$\\[1ex]
   &  $4$ &$2.6371E+04$&	$2.0138E+04$&	$3.5456E+02$	&$2.6052E+02$&	$2.7573E+04$	&$1.8454E+04$\\	[1ex]
   \midrule
   \multirow{3}{*}{$\mathrm{L}_m$}  & $2$ &$8.4606E+03$	&$8.3300E+03$&	$1.5279E+02$&	$1.1486E+02$	&$2.5276E+02$&	$1.5433E+02$\\[1ex]
   &$3$ &	$1.8511E+04$&	$1.4614E+04$	&$2.2517E+02$&	$1.6131E+02$&	$9.9621E+02$&	$5.3614E+02$\\[1ex]
   &  $4$ &$3.2837E+04$&	$2.7473E+04$	&$8.9167E+02$	&5.8061E+02&	$2.5284E+05$&	$6.7124E+04$\\	[1ex]
   \midrule
   \multirow{3}{*}{$\mathrm{L}_p$}  & $2$ &$1.1952E+04$	&$1.1786E+04$	&$2.1568E+02$	& $1.6042E+02$&	$1.0211E+03$	&$7.8013E+02$\\[1ex]
   &$3$ &$9.4332E+03$	&$7.4505E+03$	&$1.2141E+02$	&$8.3080E+01$	&$1.1802E+04$&	$8.2273E+03$	\\[1ex]
   &  $4$ &$1.8803E+04$	&$1.5754E+04$	&$7.0206E+02$	&$4.8744E+02$&	$7.4731E+03$&	$4.8946E+03$\\	[1ex]
     
   \bottomrule

    \end{tabular}
    \end{sidewaystable}

    The structured partial NCN is calculated using Theorem \ref{Th:SNCN}, while the structured partial MCN and CCN are computed from Theorem \ref{Th:MCN}. Numerical results for various choices of \(\mathrm{L}\) and \(q = 2, 3, 4\) are summarized in Table \ref{tab5}. We observed for all choices of $\mathrm{L},$ $\mathfrak{K}^{\infty,\mathbb{S}}_{mix,  \vp}(\mathbf{H}, \b; \mathrm{L})$ and $\mathfrak{K}^{\infty,\mathbb{S}}_{com,  \vp}(\mathbf{H}, \b; \mathrm{L})$  are almost one order smaller than the $\mathfrak{K}^{\infty}_{mix,  \vp}(\mathbf{H}, \b; \mathrm{L})$ and $\mathfrak{K}^{\infty}_{com,  \vp}(\mathbf{H}, \b; \mathrm{L})$, respectively. 
\end{exam}
\section{Conclusion}\label{sec:Conclusion}
This paper introduces a unified framework for investigating the partial CN for the solution of the DSPP. We derive compact formulas for the partial CNs, and by considering specific norms, we obtain the partial NCN, MCN and CCN for the solution of the DSPP. Additionally, we provide sharp upper bounds for the partial CNs that are free from costly Kronecker products. Moreover, we compute structured partial NCN, MCN, and CCN  by introducing perturbations that maintain the structure of the block matrices of the coefficient matrix. Using our theoretical findings and by leveraging the relationship between the EILS problems and the DSPP, we recover several previously established results for the EILS problems. Experimental results demonstrate that the derived upper bounds for the partial CNs provide tight estimates of the actual partial CNs. Furthermore, the proposed partial CNs and their upper bounds provide sharp error estimation for the solution, highlighting their effectiveness and reliability.

\bmhead{Acknowledgements}
During this work, Pinki Khatun was financially supported by the Council of Scientific $\&$ Industrial Research (CSIR) in New Delhi, India, in the form of a fellowship (File no. $09/1022(0098)/2020$-EMR-I).

\section*{Declarations}
\subsection*{Funding} 
Not applicable
\subsection*{Competing interests}
The authors have no competing interests
\subsection*{Ethics approval} 
Not applicable
\subsection*{Data availability} 
Not applicable
\subsection*{Materials availability} 
Not applicable
\subsection*{Author contribution}

Sk. Safique Ahmad: Contributed to the conceptualization and design of the study, provided critical feedback and suggestions on earlier versions of the manuscript, contributed to improvements, and reviewed and approved the final version of the manuscript.
Pinki Khatun: Led the conceptualization and design of the study, drafted the initial version of the manuscript, developed and tested the software, provided comments on earlier versions, contributed to revisions and improvements, and reviewed and approved the final version of the manuscript.
\bibliography{sn-bibliography}

\end{document}